\documentclass{svmult}

\usepackage{mathptmx}        % selects Times Roman as basic font
\usepackage{helvet}          % selects Helvetica as sans-serif font
\usepackage{courier}         % selects Courier as typewriter font
\usepackage{type1cm}        % activate if the above 3 fonts are 
\usepackage{makeidx}         % allows index generation
\usepackage{graphicx}        % standard LaTeX graphics tool
\usepackage{multicol}        % used for the two-column index
\usepackage[bottom]{footmisc}% places footnotes at page bottom

\usepackage{amsmath,amssymb,bbm}
\usepackage{graphicx}
\usepackage[utf8]{inputenc}
\usepackage{url}
\usepackage{microtype}
\usepackage{import}
\usepackage{color}
\usepackage{epic}
\usepackage{tikz}
\usepackage{algpseudocode}
\usepackage{setspace}
\usepackage{tabularx}
\usepackage{array}

\usepackage{amsfonts}
\usepackage{psfrag}
\usepackage{epsfig}
\usepackage{multirow}
\usepackage{bm}
\usepackage{cite}
\usepackage{float}

\usepackage{mydef}

\DeclareMathAlphabet{\mathcal}{OMS}{cmsy}{m}{n}

\let\jmath=\undefined
\DeclareSymbolFont{jmathcmletters}{OML}{cmm}{m}{it}
\DeclareMathSymbol{\jmath}{\mathord}{jmathcmletters}{"7C}

\usepackage{bm}
\usepackage{lineno} %\linenumbers*[1]
\usepackage{comment}
\usepackage{enumerate}

\newcommand{\Seqn}{\begin{equation}}
\newcommand{\Feqn}{\end{equation}}
\newtheorem{algorithm}{Algorithm}[section]

\newcommand*\samethanks[1][\value{footnote}]{\footnotemark[#1]}

\title*{A sparse grid method for Bayesian uncertainty quantification with application to large eddy simulation turbulence models\thanks{This material is based upon work supported in part by the U.S.~Air Force of Scientific Research under grant numbers 1854-V521-12; by the U.S.~Department of Energy, Office of Science, Office of Advanced Scientific Computing Research, Applied Mathematics program under contract numbers ERKJ259, ERKJE45; and by the Laboratory Directed Research and Development program at the Oak Ridge National Laboratory, which is operated by UT-Battelle, LLC., for the U.S.~Department of Energy under Contract DE-AC05-00OR22725.}
}

\titlerunning{A sparse grid method for Bayesian inference of LES turbulence models}

\authorrunning{Hoang Tran, Clayton G.~Webster, and Guannan Zhang}

\begin{document}

\author{
Hoang Tran\thanks{Computer Science and Mathematics Division, Oak Ridge National Laboratory, Oak Ridge TN 37831-6164 (\texttt{tranha@ornl.gov}, \texttt{webstercg@ornl.gov}, \texttt{zhangg@ornl.gov}).} 
\and Clayton G.~Webster\samethanks
\and Guannan Zhang\samethanks
}

%\graphicspath{{figures/}} % figures in subdirectory
%\input{author}

% We added this to correct the equation numbering within
\numberwithin{equation}{section}
\renewcommand{\theequation}{\thesection.\arabic{equation}}

\maketitle

\newcommand{\tabincell}[2]{\begin{tabular}{@{}#1@{}}#2\end{tabular}}

\abstract{
There is wide agreement that the accuracy of turbulence models suffer from their sensitivity with respect to physical input data, the uncertainties of user-elected parameters, as well as the model inadequacy. However, the application of Bayesian inference to systematically quantify the uncertainties in parameters, by means of exploring posterior probability density functions (PPDFs), has been hindered by the prohibitively daunting computational cost associated with the large number of model executions, in addition to daunting computation time per one turbulence simulation. In this effort, we {perform in this paper} an \textit{adaptive hierarchical sparse grid} surrogate modeling approach to Bayesian inference of large eddy simulation (LES). First, an adaptive hierarchical sparse grid surrogate for the output of forward models is constructed using a relatively small number of model executions. Using such surrogate, the likelihood function can be rapidly evaluated at any point in the parameter space without simulating the computationally expensive LES model. {This method is essentially similar to those developed in \cite{ZLY+13} for geophysical and groundwater models, but is adjusted and applied here for a much more challenging problem of uncertainty quantification of turbulence models. Through a numerical demonstration of the Smagorinsky model of two-dimensional flow around a cylinder at sub-critical Reynolds number, our approach is proven to significantly reduce the number of costly LES executions without losing much accuracy in the posterior probability estimation.} Here, the model parameters are calibrated against synthetic data related to the mean flow velocity and Reynolds stresses at different locations in the flow wake. {The influence of the user-elected LES parameters on the quality of output data will be discussed.  }
}

%--------------
%keywords
%

%\keyphrases{uncertainty quantification, large eddy simulation, computational fluid dynamics, eddy viscosity model}
%\AMclass{1234.56}

\section{Introduction}
For most turbulent flows encountered in industrial applications, the cost of direct numerical simulation (DNS) would exceed the capacity of current computational resource (and possibly continue to do so for the foreseeable future). As a result, many important decisions affecting our daily lives (such as climate policy, biomedical device design, pollution dispersal and energy efficiency improvement) are informed from simulations of turbulent flows by various models of turbulence. The {accuracy of estimated quantities of interest (QoIs)} by such models, however, frequently suffers from the uncertainties on the physical input data, user-chosen model parameters and the subgrid model. It is ideal to be able to incorporate these uncertainties in the predictions of QoIs. 

%Consider the Navier-Stokes equation (NSE) for the pointwise velocity $u(x,t)$ and pressure $p(x,t)$ of an incompressible, viscous, Newtonian fluid 

%\begin{gather}
%\label{NSE}
%\begin{aligned}
%&u_t + u\cdot\nabla u - \nu\Delta u + \nabla p = f,\ \mbox{ in }\Omega\times (0,T),
%\\
%&\qquad\qquad\nabla\cdot u = 0,\mbox{ in }\Omega\times (0,T),\ 
%\\
%&\quad u = 0 \mbox{ on }\partial\Omega\mbox{ and }u(x,0) = u_0(x),\mbox{ in }\Omega. 
%\end{aligned}
%\end{gather}

The basic approach used for approximating turbulent flows has been to compute the time- and space-filtered velocity and pressure, which are less computationally demanding and of main technical interest, instead of solving for the pointwise velocity and pressure prescribed by the standard Navier-Stokes equations. %In other words, turbulence models seek to simulate the large scales of a flow while suppressing any fluctuations below the cut-off length scale. Since the information about the small scales is not available, the influence of them on large scales needs to be hypothesized and modeled. Unfortunately, in many cases, the assumptions which the model bases on are poor and over-simplified. 
{The use of turbulence models leads to a level of uncertainty in the performance and inaccuracy in the simulation results, due to user-chosen model parameters whose true or optimal values are not well-known \textit{(parametric uncertainty)}, or the inherent inability of the model to reproduce reality \textit{(structural uncertainty)}.} With the fast growth in available computational power, the literature on uncertainty quantification for fluid mechanics modeling has grown extensively recently. Many stochastic numerical methods have been developed, analyzed and tested for simulations of fluid flows with uncertain physical and model parameters, see, e.g., \cite{WK06,GZ07,FWK08,WLSB08,SM11,TTW14}. Sensitivity analysis of LES to parametric uncertainty was conducted in \cite{LMS07}. Statistical methods to capture structural uncertainties in turbulence models were presented in \cite{GELI12,ELI13,GI13}. For inverse uncertainty quantification, we refer to \cite{COP+11, OM11}  (Bayesian inference for Reynolds-averaged Navier Stokes (RANS) models) and \cite{DW11} (adjoint based inverse modeling). 

Bayesian inference has become a valuable tool for estimation of parametric and structural uncertainties of physical systems constrained by differential equations. %The approach is centered around the use of Bayes' theorem and the incorporation of multiple types of data and prior information to update probability density functions of quantities of interest, which can be directly used for uncertainty quantification, risk assessment, and decision making. 
Sampling techniques, such as Markov chain Monte Carlo (MCMC), have frequently been employed in Bayesian inference \cite{RC04,GL06,LCK10}. {However, MCMC methods \cite{HLMS06,VBC+08,VBD+09} are, in general, computationally expensive, because a large number of forward model simulations is needed to estimate the PPDF and sample from it.} Given the fact that one solution of turbulence models easily takes thousands of computing hours, MCMC simulations in many CFD applications would require prohibitively large computational budgets. Perhaps due to this demand, efforts on model calibration up until now have been limited on the least expensive turbulence model $-$ RANS equations \cite{COP+11, OM11}. %It is worth noting that in the aforementioned works, conventional MCMC has been utilized to explore the posterior distribution and neither the efficiency nor accuracy of the sampling methods has been discussed. 
To make Bayesian inference tractable for other types of closure models, including LES, it is essential to perform the MCMC sampling in a time and cost effective manner. 

A strategy to improve the efficiency of MCMC simulations is \textit{surrogate modeling}, which has been developed in a wide variety of contexts and disciplines, see \cite{RTB12} and the reference therein. Surrogate modeling practice seeks to approximate the response of an original function (model outputs or the PPDF in this work), which is typically computationally expensive, by a cheaper-to-run surrogate. {The PPDF can then be evaluated by sampling the surrogate directly without forward model executions. Compared to conventional MCMC algorithms, this approach is advantageous that it significantly reduces the number of forward model executions at a desired accuracy and allows sampling the PPDF in parallel. Several methods can be employed to construct the surrogate systems, including polynomial chaos expansion \cite{GS91}, stochastic Galerkin \cite{BNZ04}, stochastic collocation \cite{BNT07}, and polynomial dimensional decomposition \cite{Rah08}, to list a few. For problems where the quantities of interest have irregular dependence with respect to the random parameters, such as those studied herein, it should be noted that approximation approaches that use global polynomials are generally less effective than those allowing for multi-level, multi-scale decomposition. In this direction, one can develop multi-level hierarchical subspaces and employ adaptive grid refinement to concentrate grid points on the subdomains with a locally high variation of solutions, resulting in a significant reduction in the number of grid points.   }

In this paper, we {present} an adaptive hierarchical sparse grid (AHSG) surrogate modeling approach to Bayesian inference of turbulence models, in particularly LES. {The key idea is to place a grid in the parameter space with sparse parameter samples, and the forward model is solved only for these samples. Compared to the regular full grid approach, sparse grid preserves the high level of accuracy with less computational work, see \cite{NTW08,NTW08b,GZ07,GWZ13b,GWZ14,BCWZ14}. As sparse grid methods require the bounded mixed derivative property, which is open for the solutions of Navier-Stokes equations and turbulence models in general, a locally adaptive refinement method, guided by hierarchical surpluses, is employed to extend sparse grid approach to possible non-smooth solutions. This refinement strategy is different from dimension-adaptive refinement \cite{GG03}, which puts more points in dimensions of higher relevance and more in line with those in \cite{Gri98,Pfl10,LGH11}. Although similar surrogate methods has been studied in \cite{MZ09, ZLY+13} for geophysical and groundwater models}, we tackle here a more challenging problem of uncertainty quantification of turbulence models. Indeed, turbulent flows are notorious for their extremely complex nature and {the non-smoothness of the surface of LES output data may weaken the accuracy of the surrogate}. The applicability of surrogate modeling techniques to LES therefore needs thorough investigation. In this work, we will demonstrate the accuracy and efficiency of the surrogate model through a numerical example of the classical Smagorinsky closure model of turbulent flow around a circular cylinder at a sub-critical Reynolds number ($Re= 500$), {which is a benchmark test case for LES}. The computation will be conducted for the two-dimensional flow, {whose outputs have similar patterns as three-dimensional simulation, but which is significantly less demanding in computing budget}. The {synthetic} data of velocity and Reynolds stresses at different locations in the flow wake are utilized for the calibration. 

This work is only one piece in the complete process of calibration and validation of LES models to issue predictions of QoIs with quantified uncertainties, and many open questions remain. We do not attempt to fit the numerical solutions with physical data herein, as the two-dimensional model has been known to show remarkable discrepancy with the experiment results. Applying our framework to the three-dimensional simulation for parameter calibration against real-world data would be the next logical step. {Another important problem is to evaluate and compare the performance of our AHSG with other surrogate methods (including some listed above) in this process. This would be conducted in future research.} {Also,} characterization and quantification of the structural inadequacy and comparison of different competing LES models are beyond the scope of this study. 

The rest of the paper is organized as follows. The Bayesian framework and the adaptive hierarchical sparse grid method of constructing the surrogate system are described in \S \ref{sec:surrogate}. In \S \ref{sec:les}, we give a detailed description of the Smagorinsky model of sub-critical flow around a cylinder. The performance of surrogate modeling approach and results of the Bayesian analysis are presented in \S \ref{sec:result}. Finally, discussions and conclusions appear in \S \ref{sec:conclusion}. 

\section{Adaptive hierarchical sparse grid methods for surrogate modeling in Bayesian inference}
%\section{Bayesian inference and adaptive sparse grid stochastic collocation method for construction of the surrogate PPDF}
\label{sec:surrogate}

\subsection{Bayesian inference}

Consider the Bayesian inference problem for a turbulence model
\begin{align}\label{def}
  \bm{d} = \bm{f}(\bm{\theta} ) + \bm{\varepsilon},
\end{align}
where $\bm{d} = (d_1,\ldots, d_{N_d})$ is a vector of $N_d$ reference data,
$\bm{\theta} = (\theta_1, \ldots,\theta_{N_\theta})$ is a vector of $N_\theta$ model parameters,
$\bm{f}(\bm{\theta})$ is the forward model, e.g., Smagorinsky model (see \S\ref{sec:les}), 
with $N_{\theta}$ inputs and $N_d$ outputs, and $\bm{\varepsilon}$ is a vector of residuals, including measurement, model parametric and structural errors. {(Nonlinear model $  \bm{d} ={\bm \Xi}(\bm{f},\bm{\theta},\bm{\varepsilon})$ can be considered as well, but leads to more complicated likelihood functions, as $\bm{\varepsilon} = \bm{\Xi}^{-1}(\bm{f},\bm{\theta})(\bm{d})$).}

The posterior distribution $P(\bm{\theta}|\bm{d})$ of the model parameters
$\bm{\theta}$, given the data $\bm{d}$, can be estimated using the Bayes' theorem \cite{BT92} via
\begin{linenomath*}\begin{equation}\label{s2:e2}
  P(\bm{\theta}|\bm{d}) = \frac {L(\bm{\theta}|\bm{d})P(\bm{\theta})}{\int L(\bm{\theta}|\bm{d})P(\bm{\theta})d\bm{\theta}},
\end{equation}\end{linenomath*}
where $P(\bm{\theta})$ is the prior distribution and $L(\bm{\theta}|\bm{d})$ is the likelihood function that measure ``goodness-of-fit'' between model simulations and observations. In parametric uncertainty quantification, the denominator of the Bayes' formula in equation \eqref{s2:e2} is a normalization constant that does not affect the shape of the PPDF. As such, in the hereafter discussion concerning building surrogate systems,
the notation $P(\bm{\theta}|\bm{d})$ or the terminology PPDF will only refer to the product $L(\bm{\theta}|\bm{d})P(\bm{\theta})$. The prior distribution represents knowledge of the parameter values before the data $\bm{d}$ is available. {When prior information is lacking, a common practice is to assume uniform distributions with parameter ranges large enough to contain all plausible values of parameters.}

Selection of appropriate likelihood functions for a specific turbulence simulation is an open question. A commonly used \textit{formal} likelihood function is based on the {simplistic} assumption that the residual term $\bm{\varepsilon}$ in \eqref{def} follows a multivariate Gaussian distribution with mean zero and prescribed standard deviations, which leads to the Gaussian likelihood function: 
\begin{align}
L(\bm{\theta}|\bm{d}) = \exp\left[-\frac{1}{2}(\bm{d}-\bm{f}(\bm{\theta}))^{\top} \Sigma^{-1}(\bm{d}-\bm{f}(\bm{\theta})) \right]. \label{MVN} \tag{MVN}
\end{align}
In this paper, we assume that the residual errors are independent, i.e., the covariance matrix $\Sigma$ is diagonal. {To describe the correlation of the errors or the inadequacy of turbulence models, other covariance matrices can also be used (and lead to inconsistent results) \cite{COP+11, OM11}.} In general, the formal approach has been criticized for relying heavily on residual error assumptions that do not hold. Alternatively, \textit{informal} likelihood functions are proposed as a pragmatic approach to implicitly account for errors in measurements, model inputs and model structure and to avoid over-fitting to reference data \cite{BB92}. Definition of informal likelihood functions is problem specific in nature, and there has been no consensus on which informal likelihood functions outperforms others. For the sake of illustration, in \S\ref{sec:result}, the exponential informal likelihood function is used for the numerical example (together with \eqref{MVN}). It reads: 
\begin{align}
L(\bm{\theta}|\bm{d}) = \exp\left(-\zeta \cdot \dfrac{\sum_{i=1}^{N_d}  \left((d_i-f_i) - (\overline{\bm{d}}-\overline{\bm{f}})\right)^2}
{\sum_{i=1}^{N_d}  \left(d_i- \overline{\bm{d}}\right)^2}\right), \label{EXP} \tag{EXP}
\end{align}
where $\overline{\bm{d}}$ is the mean of observations, $\overline{\bm{f}}$ is the mean of the outputs of forward model, and $\zeta$ is a scaling constant. For some other widely used informal likelihood functions, see \cite{SBT08}. 

\subsection{Adaptive hierarchical sparse grid methods for construction of the surrogate PPDF}

The central task of Bayesian inference is to estimate the posterior distribution $P(\bm{\theta}|\bm{d})$. It is often difficult to draw samples from the PPDF directly, so the MCMC methods, such as the Metropolis-Hastings (M-H) algorithm \cite{GL06} and its variants, are normally used for the sampling process. In practice, the convergence of MCMC methods is often slow, leading to a large number of model simulations. To tackle this challenge, surrogate modeling approaches seek to build an approximation (called the surrogate system) for $P(\bm{\theta}|\bm{d})$, then the MCMC algorithm draws samples from it directly without executing the forward model. With this approach, the main computational cost for evaluating the PPDF is now transferred to the surrogate construction step. Naturally, an approximation method which requires minimal number of grid points in the parameter space, while not surrendering much accuracy is desired. The methodology we utilize to construct the surrogate system, presented in this subsection, is similar to the method introduced in \cite{ZLY+13}. {Since the method can be applied to functions governed by partial differential equations, not limited to $P(\bm{\theta}|\bm{d})$ or $\bm{f}(\bm{\theta})$, a generic notation $\eta(\bm{\theta}):\Omega\to \mathbb{R}$ is used for the description. The following assumptions are needed: 

\begin{enumerate}[\ \ \ \ \ (a)]
\item The domain $\Omega$ is a rectangle, i.e., $\Omega = \Omega_1 \times \ldots \times \Omega_{N_\theta}, $ where $\Omega_n \subset \mathbb{R},\, n= 1,\ldots, N_\theta.$
\item The joint probability density function ${ \rho}({\bm \theta})$ is of product-type: 
\begin{align*}
{ \rho}({\bm \theta}) = \prod_{n=1}^{N_{\theta}} \rho_n(\theta_n),
\end{align*}
where $\rho_n:\Omega_n\to \mathbb{R}$ are univariate density functions. 
\item The univariate domains and density functions are identical: 
\begin{align*}
\Omega_1 = \ldots = \Omega_{N_{ \theta}};\ \rho_1 = \ldots = \rho_{N_{ \theta}},
\end{align*}
yielding the same $i$-level univariate quadrature rules 
\begin{align*}
\mathcal{Q}_i^{(1)}[\cdot] = \ldots = \mathcal{Q}_i^{(N_\theta)}[\cdot] =:  \mathcal{Q}_i[\cdot] . 
\end{align*}
\item The univariate quadrature rules are nested.
\end{enumerate}
}

\subsubsection{Adaptive sparse grid interpolation}\label{s_1dhintp}

%\begin{comment}
%
The basis of constructing the sparse grid approximation in the multi-dimensional setting is the one-dimensional (1-D) hierarchical interpolation.
Consider a function
$\eta(\theta): [0,1] \rightarrow \mathbb{R}$. The 1-D hierarchical Lagrange interpolation formula
is defined by
\begin{linenomath*}\begin{equation}\label{1d_hintp}
\mathcal{U}_K[\eta](\theta) := \sum_{i=0}^K \Delta \mathcal{U}_i[\eta](\theta),
\end{equation}\end{linenomath*}
where $K$ is the resolution level, and the incremental interpolation operator $\Delta \mathcal{U}_i[\eta]$ is given as
\begin{linenomath*}\begin{equation}\label{delta_intp}
 \Delta \mathcal{U}_i[\eta](\theta) := \sum_{j=1}^{m_i} c_{i,j} \phi_{i,j}(\theta), \quad i = 0,\ldots,K.
\end{equation}\end{linenomath*}
%
%The nonnegative integer $L$ in equation \eqref{1d_hintp} is calthe \textit{resolution level} of the
%hierarchical interpolant $\mathcal{U}_L[\eta]$ and the summation over the resolution level in
%\eqref{1d_hintp} exhibits the hierarchical structure of the interpolant $\mathcal{U}^L(\eta)$.
For $j = 1,\ldots, m_i$, $\phi_j^i(\theta)$ and $c_{i,j}$ in \eqref{delta_intp} are the piecewise hierarchical basis functions  \cite{BG04,ZLY+13}
and the interpolation coefficients for $\Delta \mathcal{U}_i[\eta]$, respectively.
For $i = 0,\ldots,K$, the integer $m_i$ in \eqref{delta_intp} is the number of interpolation points
involved in $\Delta \mathcal{U}_i[\eta]$, which is defined by
\begin{equation*}
    m_0 = 1,\;\; m_1 = 2,\;\; \text{ and }\;\;  m_i = 2^{i-1} \;\; \text{for} \;\; i \ge 2.
 \end{equation*}
A uniform grid, denoted by $\Delta \mathcal{X}_i = \{\theta_{i,j}\}_{j=1}^{m_i}$, can be utilized
for the incremental interpolant $\Delta \mathcal{U}_i[\eta]$. The abscissas of $\Delta \mathcal{X}_i$ are defined by
\[
\theta_{0,1} = 0.5,\;\;  \theta_{1,1} = 0,\;\; \theta_{1,2} =1, \;\; \text{ and }\;\; \theta_{i,j} =  \dfrac{2j-1}{\sum_{k=0}^{i}m_k-1} \;\; \text{for} \;\; j = 1\ldots,m_i, \; i \ge 2.
\]
%\begin{linenomath*}\begin{equation*}
%  \begin{array}{ll}
%     \theta_1^0 = 0.5 &\quad \text{for} \quad i = 0,\\
%     \\
%     \theta_1^1 = 0, \; \theta_2^1 = 1 &\quad \text{for} \quad i = 1,\\
%     \\
%     \theta_j^i =  \dfrac{2j-1}{{\displaystyle\sum_{k=0}^{i}}m_k-1} &\quad \text{for} \quad j = 1\ldots,m_i, \; i \ge 2.
%  \end{array}
%\end{equation*}\end{linenomath*}
%
Then, the  hierarchical grid for $\mathcal{U}_K[\eta](\theta)$ is defined by
%
%\begin{linenomath*}\begin{equation*}
$\mathcal{X}_K = \cup_{i=0}^K \Delta \mathcal{X}_i.$
%\end{equation*}\end{linenomath*}
%

%The coefficient $c_{i,j}$ is defined as the hierarchical \textit{surplus} of the basis function
%$\phi_j^i(\theta)$, which is the difference between the value of the interpolated function $\eta(\theta)$
%and the value of the interpolant $\mathcal{U}^{i-1}(\eta)$ at $\theta_j^i$. For $i = 0$,
%%
%\begin{linenomath*}\begin{equation*}
%c_1^0 = \Delta \mathcal{U}^0(\eta)(\theta_1^0) = \mathcal{U}^0(\eta)(\theta_1^0) = \eta(\theta_1^0),
%\end{equation*}\end{linenomath*}
%%
%and for $i >  0$, $j = 1,\ldots, m_i$,
%%
%\begin{linenomath*}\begin{equation*}
%\begin{aligned}
%c_j^i &= \Delta \mathcal{U}^i(\eta)(\theta_j^i) \\
%      &= \mathcal{U}^i(\eta)(\theta_j^i) - \mathcal{U}^{i-1}(\eta)(\theta_j^i)\\
%      &= \eta(\theta_j^i) - \mathcal{U}^{i-1}(\eta)(\theta_j^i).
%\end{aligned}
%\end{equation*}\end{linenomath*}
%
%\end{comment}

Based on the one-dimensional hierarchical interpolation, we can construct an approximation for a multivariate function $\eta(\bm{\theta}):[0,1]^{N_{\theta}} \rightarrow \mathbb{R}$, where $\bm{\theta} =
(\theta_1,\ldots, \theta_{N_\theta})$, by hierarchical interpolation formula as
\begin{linenomath*}\begin{equation}\label{md_hintp}
\mathcal{I}_{K}[\eta](\bm{\theta}) := \sum_{|{\bf i}| \le K} \Delta_{{\bf i}}[\eta](\bm{\theta})
\end{equation}\end{linenomath*}
and the multi-dimensional incremental interpolation operator $\Delta_{{\bf i}}[\eta]$ is defined by
\begin{linenomath*}\begin{equation*}%\label{md_delta_intp}
\begin{aligned}
\Delta_{{\bf i}}[\eta](\bm{\theta}) & := \Delta \mathcal{U}_{i_1}
\otimes \cdots \otimes \Delta \mathcal{U}_{i_{N_\theta}} [\eta](\bm{\theta})= \sum_{{\bf j}\in B_{{\bf i}}}c_{\ii,\jj} \bm{\phi}_{\ii,\jj}(\bm{\theta}),%\\
%& = \sum_{j_1=1}^{m_{i_1}} \cdots \sum_{j_{N_\theta}=1}^{m_{i_{N_\theta}}} c_{{\bf j}}^{{\bf i}}
%\bm{\phi}_{{\bf j}}^{{\bf i}}(\bm{\theta})
\end{aligned}
\end{equation*}\end{linenomath*}
where ${\bf i} := (i_1,\ldots, i_{N_\theta})$ is a multi-index indicating the resolution level of
$\Delta_{{\bf i}}[\eta]$, $ |{\bf i}| = i_1 + \cdots + i_{N_\theta}$, 
$\bm{\phi}_{\ii,\jj}(\bm{\theta}) := \prod_{n=1}^{N_\theta}\phi_{i_n,j_n}(\theta_n)$, and the multi-index set $B_{{\bf i}}$ is defined by 
$
B_{{\bf i}}= \left.\left\{{\bf j} \in \mathbb{N}^{N_\theta} \right| j_n = 1,\ldots,m_{i_n}, n = 1,\ldots,N_\theta \right\}
$.
As such, the grids for $\Delta_{\ii}[\eta]$ and $\mathcal{I}_K[\eta]$ are defined by $\Delta \mathcal{H}_{\ii}:= \Delta \mathcal{X}_{i_1} \times \cdots \times  \Delta \mathcal{X}_{i_{N_\theta}}$ and $\mathcal{H}_{K}:=  \cup_{|\ii| \le K} \Delta \mathcal{H}_{\ii}$.

{In this paper, we employ the piecewise linear hierarchical basis \cite{BG04,ZLY+13}} and the surplus $c_{\ii,\jj}$ can be explicitly computed as
\begin{linenomath*}\begin{align*}
&c_{\bm{0},\bm{1}}  = \Delta_{\bm{0}}[\eta](\bm{\theta}_{\bm{0},\bm{1}})
 = \mathcal{I}_{0}[\eta](\bm{\theta}_{\bm{0,1}}) = \eta(\bm{\theta}_{\bm{0,1}}), \\
 &c_{\ii,\jj}  =  \Delta_{\bf i}[\eta](\bm{\theta}_{\ii,\jj})
            %  = \mathcal{I}^{L,N_\theta}(\eta)(\bm{\theta_{{\bf j}}^{{\bf i}}})
             %     - \mathcal{I}^{L-1,N_\theta}(\eta)(\bm{\theta_{{\bf j}}^{{\bf i}}})
               = \eta(\bm{\theta}_{\ii,\jj}) - \mathcal{I}_{K-1}[\eta](\bm{\theta}_{\ii,\jj})\;\;\mbox{ for $|{\bf i}| = K >0$},
\end{align*}\end{linenomath*}
{as the supports of basis functions are mutually disjoint on each subspace.} As discussed in \cite{BG04}, when the function $\eta(\bm{\theta})$ is smooth with respect to $\bm{\theta}$, the
magnitude of the surplus $c_{\ii,\jj}$ will approach to zero as the resolution level $K$ increases. Therefore, the surplus
can be used as an error indicator for the interpolant $\mathcal{I}_{K}[\eta]$ in order to detect the smoothness
of the target function and guide
the sparse grid refinement. In particular, each point $\bm{\theta}_{\ii,\jj}$ of the isotropic level-$K$ sparse grid $\mathcal{H}_{K}$ is assigned two children in each $n$-th direction, represented by
\begin{linenomath*}\begin{equation}
\begin{aligned}
&C_1^n(\bm{\theta}_{\ii,\jj}) = \left(\theta_{i_1,j_1}, \ldots, \theta_{i_{n-1}, j_{n-1}},
\theta_{i_n+1, 2j_n-1}, \theta_{i_{n+1}, j_{n+1}}, \ldots, \theta_{i_{N_\theta}, j_{N_\theta}}\right),\\
&C_2^n(\bm{\theta}_{\ii,\jj}) = \left(\theta_{i_1,j_1}, \ldots, \theta_{i_{n-1},j_{n-1}},
\theta_{i_n+1,2j_n}, \theta_{i_{n+1},j_{n+1}}, \ldots, \theta_{i_{N_\theta},j_{N_\theta}}\right),\\
\end{aligned}
\end{equation}\end{linenomath*}
for $n = 1, \ldots, N_\theta$. Note that %$C_1^n(\bm{\theta}^{{\bf i}}_{{\bf j}}), C_2^n(\bm{\theta}^{{\bf i}}_{{\bf j}}) \in \Delta \mathcal{H}^{\widehat{{\bf i}}_n, N_\theta}$ with $\widehat{{\bf i}}_n = (i_1,\ldots, i_{n-1}, i_n+1, i_{n+1}, \ldots,i_{N_\theta})$ and $|\widehat{{\bf i}}_n| = |{\bf i}|+1$, i.e., 
the children of each sparse grid point on level $|{\bf i}|$ belong to
the sparse grid point set of level $|{\bf i}|+1$. The basic idea of adaptivity is as follows: for each point whose
magnitude of the surplus is larger than the prescribed error tolerance, we refine the grid by adding its children on the next level. More rigorously, for an error tolerance $\alpha$, the adaptive sparse grid interpolant is defined on each successive interpolation level as
\begin{linenomath*}\begin{equation}\label{md_aintp}
\mathcal{I}_{K,\alpha}[\eta](\bm{\theta}) := \sum_{|{\bf i}|\le K} \sum_{{\bf j} \in B_{{\bf i}}^{\alpha}}
c_{\ii,\jj} \bm{\phi}_{\ii,\jj}(\bm{\theta}),
\end{equation}\end{linenomath*}
where the multi-index set $B_{{\bf i}}^\alpha$ is defined by modifying the multi-index set $B_{{\bf i}}$, i.e.,
%
%\begin{linenomath*}\begin{equation*}%\label{ind_set2}
$B_{{\bf i}}^\alpha = \{{\bf j} \in B_{{\bf i}} | |c_{\bf j}^{\bf i}| > \alpha \}$.
%\end{equation*}\end{linenomath*}
%
The corresponding adaptive sparse grid is a sub-grid of the level-$K$ isotropic sparse grid $\mathcal{H}_{K}$, with the grid points becoming concentrated in the non-smooth region. In the region where
$\eta(\bm{\theta})$ is very smooth, this approach saves a significant number
of grid points but still achieves the prescribed accuracy. 
 
\subsubsection{Algorithm for constructing the surrogate PPDF}

In the forthcoming numerical illustration, a surrogate PPDF will be constructed based on the sparse grid method, discussed above, 
with the use of the following procedure. 

\begin{algorithm}$ $
\label{surrogate_construct}
\begin{itemize}

\item STEP 1: Determine the maximum allowable resolution $K$ of the sparse grid by analyzing the trade off between the interpolation error and computational cost. Determine the error tolerance $\alpha$. 

\item STEP 2: Generate the isotropic sparse grid at some starting coarse level $\ell$. Until the maximum level $K$ is reached or the magnitudes of all surpluses on the last level are smaller than $\alpha$, do the following iteratively: 

\begin{itemize}
\item Step 2.1: Simulate the turbulence model $\bm{f}(\bm{\theta})$ at each grid point $\bm{\theta}_{{\bf i}, {\bf j}} \in \mathcal{H}_{\ell}$. 
\item Step 2.2: Construct the sparse grid interpolant $\mathcal{I}_{\ell,\alpha}[\bm{f}](\bm{\theta}) $ based on formula \eqref{md_aintp}.
\item Step 2.3: Generate the adaptive sparse grid for the next level based on the obtained surpluses. Set $\ell:=\ell +1$ and go back to Step 2.1.
\end{itemize}  

\item STEP 3: Construct an approximate likelihood function, denoted by $\tilde{L}(\bm{\theta}|\bm{d})$, by substituting $\mathcal{I}_{\ell, \alpha}[\bm{f}] $ for $\bm{f}$ into the likelihood formula using, e.g., \eqref{MVN} or \eqref{EXP}.
\item STEP 4: Construct the surrogate PPDF $\tilde{P}(\bm{\theta}|\bm{d})$ via 
\begin{align*}
 \tilde{P}(\bm{\theta}|\bm{d}) \propto \tilde{L}(\bm{\theta}|\bm{d})P(\bm{\theta}). 
\end{align*}

\end{itemize}
\end{algorithm}

After the surrogate is constructed, an MCMC simulation is used to explore $\tilde{P}(\bm{\theta}|\bm{d})$. Using our approach, drawing the parameter samples does not require any model executions but negligible computational time for polynomial evaluation using the surrogate system. The improvement of computational efficiency by using surrogate PPDF is more impressive when increased samples are drawn in the MCMC simulation. 

Finally, it is worth discussing the flexibility of grid adaptive refinement strategies. It is known that in calibration problems of turbulence models, different likelihood models could lead to conflicting posterior distributions \cite{COP+11,OM11}. Moreover, for a flow problem, experimental data given by different authors is sometimes inconsistent. There is also a wide variation of the physical quantities to be measured and recorded. Naturally, one would desire a surrogate modeling method that allows for the use of a variety of 
likelihood functions and data sets, at little cost, once the surrogate system has been built. An adaptive refinement strategy based on the smoothness of the likelihood functions \cite{ZLY+13} is obviously the least flexible, since the grid is likelihood-function-specific. The approach we apply in this work, i.e., an adaptive method that is guided by the smoothness of output interpolant, 
allows the use of an universal surrogate of the output, for different choices of likelihood functions and data of the \textit{same} physical quantities. The surrogate for the output is, however, more expensive than that built directly for the likelihood function in the former approach, since grid points may be generated in the low density region of the likelihood where the forward simulations are wasteful. The most versatile method is certainly the non-adaptive, full sparse grid method, but the surrogate is also constructed with highest cost in this case. To this end, one has to sacrifice the flexibility of the sparse grid surrogate to improve the efficiency. The demand of investigating posterior distribution over different likelihood functions and data sets and the computational budget need to be balanced before an adaptive refinement strategy is determined. 

\section{Application to large eddy simulation of sub-critical flow around a circular cylinder}
\label{sec:les}

\subsection{{Parametric uncertainty of Smagorinsky model}}

%It is widely accepted that the Navier-Stokes equation could describe even the most complex turbulent flow. However, due to a large range of scales present, a DNS of most engineering flow problems is infeasible for the foreseeable future. This fact, coupled with the observation that for industrial applications one is usually concerned with only the large scale behavior of the flow, leads to consideration of less complete yet cheaper turbulence models. 

%The concept of LES is motivated by the idea that the physical dynamics of the large and small eddies of turbulent flows are different: the large structures evolve deterministically and are not sensitive, while the small eddies are sensitive but they have universal features so that their mean effects on the large eddies can be modeled. LES seeks to directly compute the large turbulent scales and only model the small scales. 
{
In LES practice, the time dependent, incompressible Navier-Stokes equations are filtered by, e.g., box filter, Gaussian filter, differential filter and the governing equations are given by
\begin{gather}
\label{LES-EV}
\begin{aligned}
&\overline{u}_t + \nabla\cdot(\overline{u}\ \overline{u}) - \nu\Delta \overline{u} + \nabla \overline{p} -\nabla\cdot(2\nu_T\nabla^s \overline{u})&= \overline{f},%\ \mbox{ in }\Omega\times (0,T),
\\
&\qquad\qquad\qquad\qquad\nabla\cdot \overline{u} &= 0,%\mbox{ in }\Omega\times (0,T),\ 
\end{aligned}
\end{gather}
where $\overline{u}$ is the velocity at the resolved scales, $\overline{p}$ is the corresponding pressure, $\nu_T\ge 0$ is the \textit{eddy viscosity} and $\nabla^s$ is the symmetric part of $\nabla$ operator, see \cite{BIL04}. }

The most common choice for $\nu_T$, which is studied herein, is known in LES as the Smagorinsky model \cite{NR50,Sma63} in which 
\begin{align}
\nu_T = \ell_S^2 |\nabla^s \overline{u}|, \label{eddyViscosity}
\end{align}
where $\ell_S = C_S\delta $ and $ |\cdot| = \sqrt{2(\cdot)_{ij}(\cdot)_{ij}}$, $\ell_S$ is called the Smagorinsky lengthscale. There are two model calibration parameters in this term - the Smagorinsky constant $C_S$ and the filter width $\delta$. The pioneering analysis of Lilly \cite{Lilly67}, under some optimistic assumptions, proposed that $C_S$ has a universal value $0.17$ and is not a ``tuning" constant. This universal value has been found later not the best choice for most LES computations and various different values ranging from $0.1$ to $0.25$ are usually selected leading to improved results, see, e.g., \cite{Dea70,CFR79, MMRF79, Ant81, BR81, MK82, Mas89, MD90}. The optimal choice for $C_S$ depends on the flow problems considered and even may be different for different regions in a flow field. Indeed, this poses a major drawback of the Smagorinsky model. 

The second calibration parameter - the filter width $\delta$ - characterizes the short lengthscale fluctuations to be removed from the flow fields. Ideally, the filter width should be put at the smallest persistent, energetically significant scale (the flow microscale), which demarcates the deterministic large eddies and isotropic small eddies, \cite{Pope00}. Unfortunately, such a choice is infeasible, since the flow microscale is seldom estimated. Instead, due to the fact that LES requires the spatial resolution $h$ to be proportional to $\delta$, the usual practice is to specify the grid to be used in the computation, and then take the filter width according to the grid size. The specification of grid and filter without knowledge of the microscale could lead to poor simulation. 

An additional calibration parameter involves in near wall treatment. The correct behavior of Smagorinsky eddy viscosity $\nu_T$ near the wall is $\nu_T \simeq 0$, since there is no turbulent fluctuation there. In contrast, the formulation \eqref{eddyViscosity} is nonzero and introduces large amounts of dissipation in the boundary layer. One approach to overcome this deficiency is to damp $\ell_S$ as the boundary is approached by the van Driest damping function \cite{vD56}. The van Driest scaling reads: 
\begin{align}
\label{vanDriest}
\ell_S = C_S\delta \left(1-e^{-{y^+}^n/{A^+}^n}\right)^p,
\end{align}
where $y^+$ is the distance from wall in wall units, $A^+$ is van Driest constant ascribed the value $A^+ = 25$. Various different values of $(n,p)$ have been used - the most commonly chosen are $(1,1)$ and $(3,0.5)$, \cite{RFBP97}. For simplicity, in this work, we fix $n=1$ and treat $p$ only as a calibration parameter. The variation of $p$ alone can capture the full spectrum of near wall scaling: $p=0$ means no damping function is applied, while a large $p$ associates with fast damping. We call $p$ van Driest damping parameter.  

\subsection{Sub-critical flow around a circular cylinder}

The flow concerned in this study corresponds to a time-dependent flow through a channel around a cylinder. External flows past objects have been the subject of numerous theoretical, experimental and numerical investigations because of their many practical applications, see \cite{BW72, Nor87, BM94} and the reference therein. In the sub-critical Reynolds number range ($300 < Re < 2 \times 10^5$), these flows are characterized by turbulent vortex streets and transitioning free shear layers.

We consider the two-dimensional flow around a cylinder 
%centered at $(0.5, 0.7)$ 
of diameter $D = 0.1$ in rectangular domain of size $2.2\times 1.4$, consisting a $5D$ upstream, $17D$ downstream and $7D$ in lateral directions. We employ the finite element method with second order Taylor-Hood finite element and polygonal boundary approximation. Our computation is carried out on triangular meshes generated based on Delaunay-Voronoi algorithm and refined around the cylinder. The ratio of number of mesh points on the top/bottom boundaries, left/right boundaries and cylinder boundary is fixed at 3:2:4. As common practice, the filter width is chosen locally at each triangle as the size of the current triangle. Its value therefore varies throughout the domain, and is roughly 10 times smaller near the cylinder than that in the far field. Since the synthetic data will be taken in the near wake region, for simplicity, we characterize $\delta$ by the value of the filter width on the cylinder surface. 
%\vspace{-0.1in}
\begin{figure}[h]
\centering
\includegraphics[scale = 0.45]{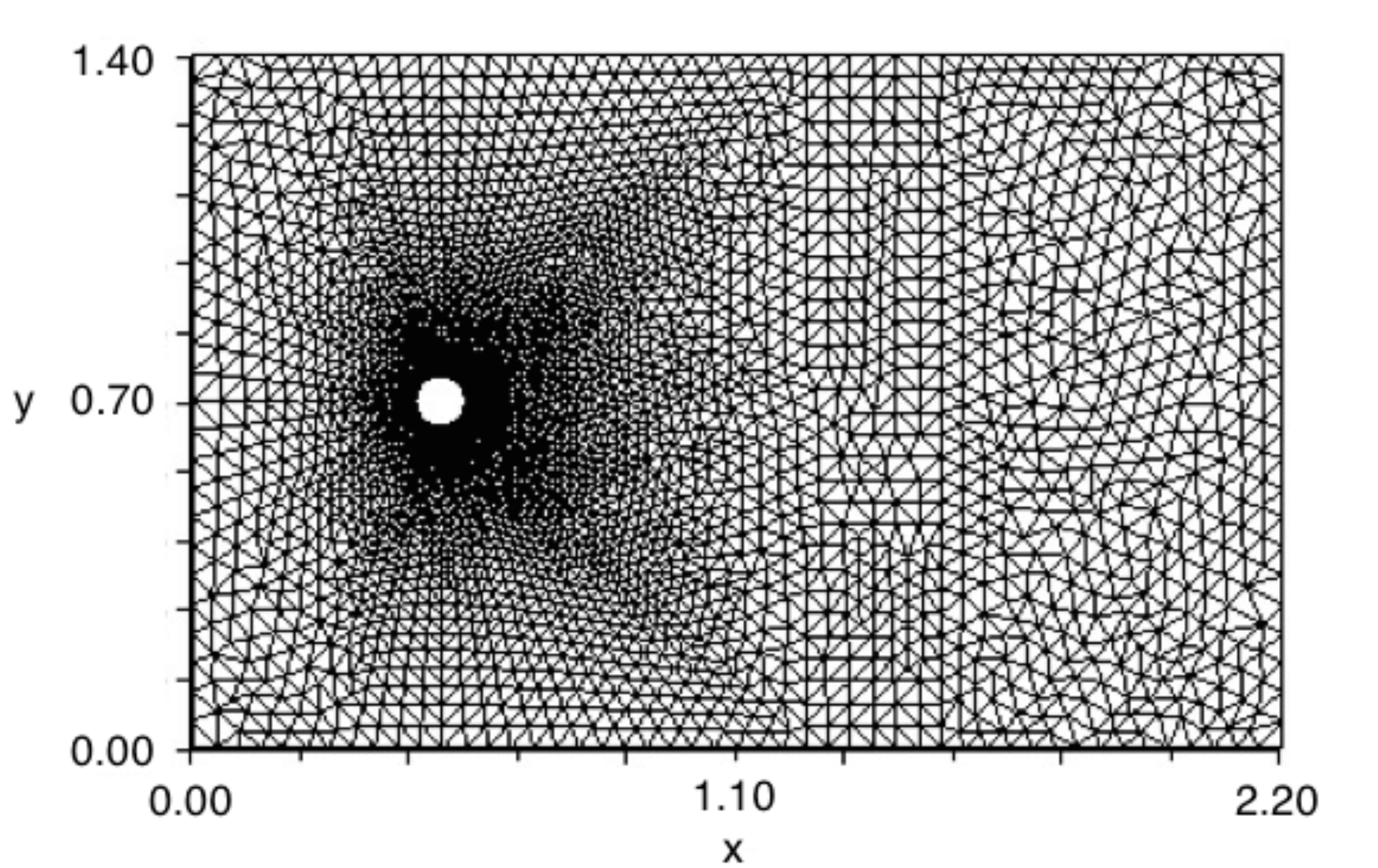}
\caption{A computational grid used in our study on LES of turbulent flow past a cylinder with $\delta = \pi/480$. }
\label{mesh}
\end{figure}
%\vspace{-0.2in}

The Smagorinsky model with van Driest damping \eqref{LES-EV}--\eqref{vanDriest} is considered with $\nu = 2\times10^{-4}$, $f = 0$, $T = 12$ and $\Delta t = 0.01$. The statistics are compiled over the last $7$ time units, equivalent to a period of $\approx 15$ vortex shedding cycles. The inflow and outflow velocity is $(\frac{6}{1.4^2} y(1.4 - y), 0)$. The mean velocity at the inlet is $U_0 =1$. No-slip boundary conditions are prescribed along the top and bottom walls. Based on $U_0$ and the diameter of the cylinder $D$, the Reynolds number for this flow is $Re = 500$, in the sub-critical range. The temporal discretization applied in the computation is the Crank-Nicolson scheme. Denoting quantities at time level $t_k$ by a subscript $k$, the time stepping scheme has the form:
\begin{gather}
\label{CrankNicolson}
\begin{aligned}
&\frac{\overline{u}_k - \overline{u}_{k-1}}{\Delta t} - \nu\Delta \frac{\overline{u}_k + \overline{u}_{k-1}}{2} + \frac{1}{2}(\overline{u}_k\cdot\nabla \overline{u}_k + \overline{u}_{k-1}\cdot\nabla \overline{u}_{k-1} ) +  \nabla \overline{p}_k 
\\
& \qquad - (\nabla\cdot(\nu_T(\overline{u}_k)\nabla^s \overline{u}_k) + \nabla\cdot(\nu_T(\overline{u}_{k-1})\nabla^s \overline{u}_{k-1})) = 0,
\\
& \nabla\cdot \overline{u}_k = 0. 
\end{aligned}
\end{gather}
System \eqref{CrankNicolson} is reformulated as a nonlinear variational problem in time step $t_k$. This problem is solved iteratively by a fixed point iteration. Let $(\overline{u}_k^0,\overline{p}_k^0)$ be an initial guess. Given $(\overline{u}^m_k,\overline{p}^m_k)$, the iterate $(\overline{u}^{m+1}_k,\overline{p}^{m+1}_k)$ is computed by solving 
\begin{gather}
\label{FixedPoint}
\begin{aligned}
&\frac{\overline{u}^{m+1}_k - \overline{u}_{k-1}}{\Delta t} - \nu\Delta \frac{\overline{u}^{m+1}_k + \overline{u}_{k-1}}{2} + \frac{1}{2}(\overline{u}_k^{m}\cdot\nabla \overline{u}_k^{m+1} + \overline{u}_{k-1}\cdot\nabla \overline{u}_{k-1} ) +  \nabla \overline{p}^{m+1}_k 
\\
& \qquad - (\nabla\cdot(\nu_T(\overline{u}_k^m)\nabla^s \overline{u}_k^{m+1}) + \nabla \cdot(\nu_T(\overline{u}_{k-1})\nabla^s \overline{u}_{k-1})) = 0,
\\
& \nabla\cdot \overline{u}^{m+1}_k = 0. 
\end{aligned}
\end{gather}
The fixed point iteration in each time step is stopped if the Euclidean norm of the residual vector is less than $10^{-10}$. The spatial and temporal discretizations we use herein are similar to \cite{JM01,John04}, in which they were applied to direct numerical simulations of flow around a cylinder at Reynolds number $Re = 100$. 

\subsection{The prior PDF and calibration data}

\begin{figure}[h]
\centering
\includegraphics[width = 2.2in]{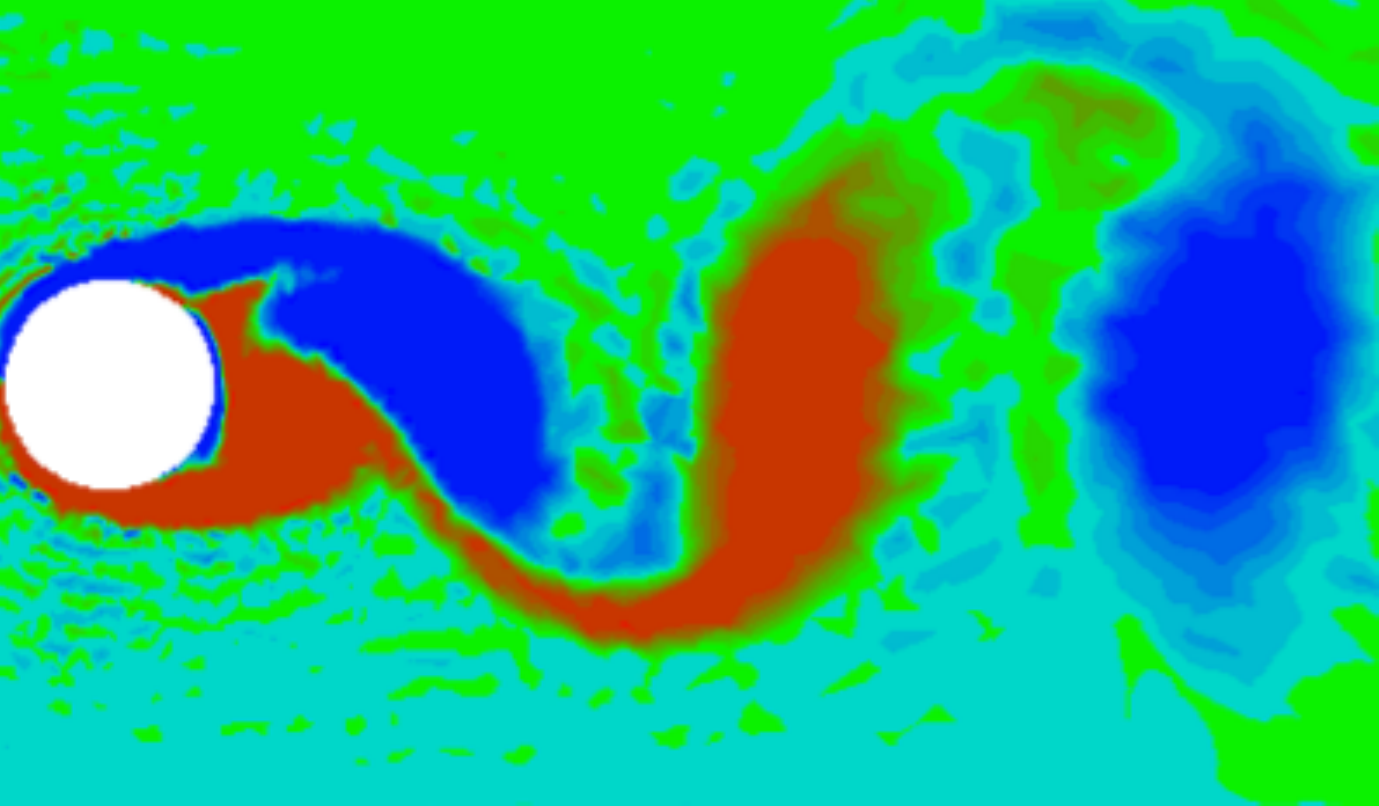}
\hspace{0.1in}
\includegraphics[width = 2.2in]{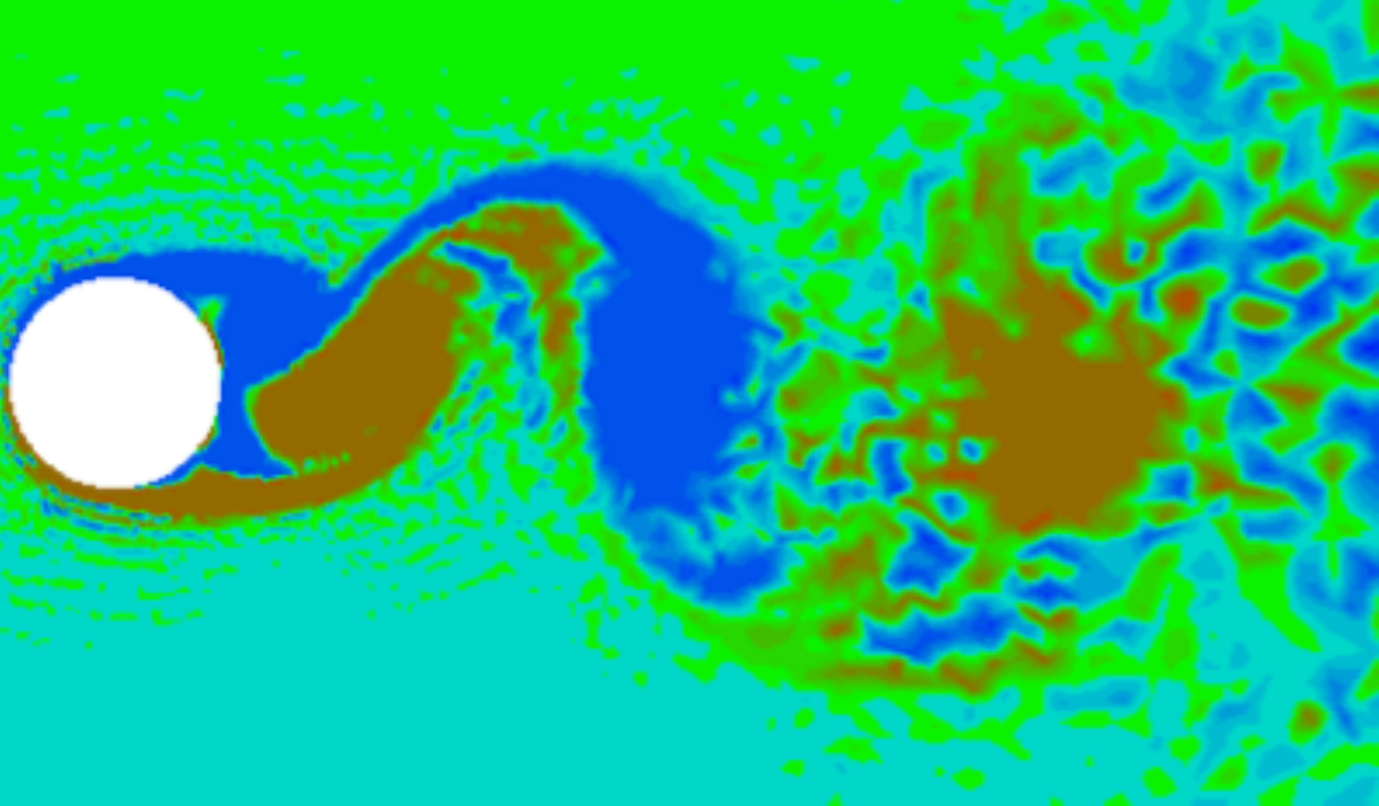}
\caption{Intantaneous vorticity at $t=20$ generated by two different choices of model parameters. Left: $C_s = 0.2,\, p = 0,\, \delta = \pi/480$. Right: $C_s = 0.05,\, p = 0,\, \delta = \pi/720$.}
%\vspace{-0.2in}
\label{inst_flow}
\end{figure}

We will exploit Bayesian calibration for three model parameters $C_S$, $p$ and $\delta$. The uniform prior PDF of the uncertain parameters is assumed. The searching domains for $C_S$ and $p$ are $[0,0.2]$ and $[0,2]$ respectively, covering their plausible and commonly selected values. The range of the prior PDF of $\delta$, on the other hand, would significantly affect the computational cost; since the filter width is proportional to the spatial resolution. Thus, to reduce the cost of flow simulations, the searching domain for $\delta$ is set to be $[\pi/600,\pi/200]$, corresponding to relatively coarse resolutions where the grid spacing on the cylinder surface ranges from $\approx 2$ to $6$ wall units. As we shall see, the response surfaces tend to be more complicated for the low-resolution simulation, possibly due to the non-physical oscillations in the underresolved solutions reflecting in the probability space. As a result, coarse grids pose a greater challenge for the surrogates to precisely describe the true outputs and are suitable for our purpose of verifying the accuracy of the surrogate modeling approach. Figure \ref{inst_flow} shows the distribution of instantaneous vorticity at $t=20$ in the near wake region for two different choices of turbulence parameters. We can see that the simulated flows display laminar vortex shedding, as expected for LES of flows past bluff bodies. The difference in phase of vortex shedding in two simulations is recognizable.       

\begin{figure}[h]
\centering
\includegraphics[height = 1.12in]{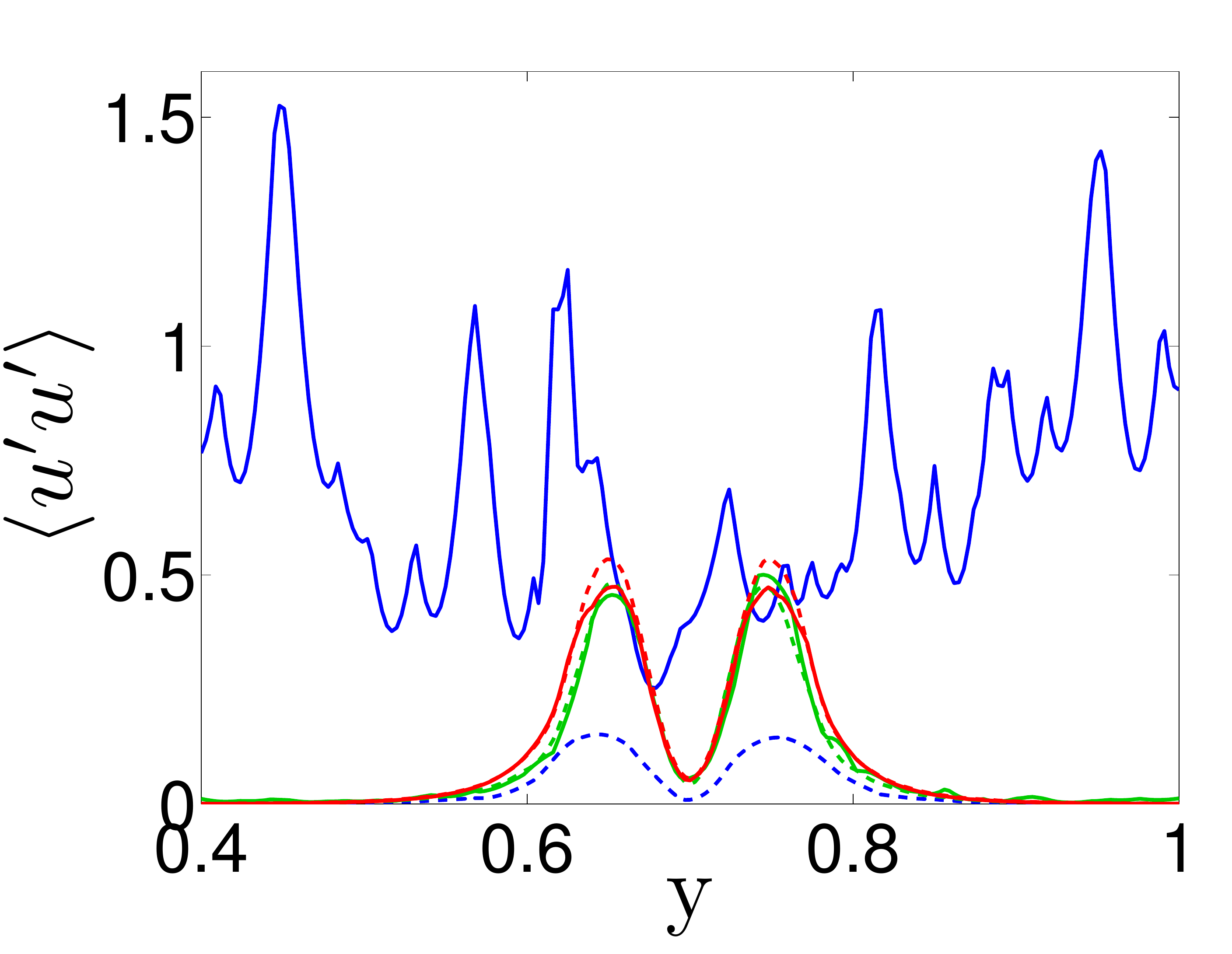}
\includegraphics[height = 1.12in]{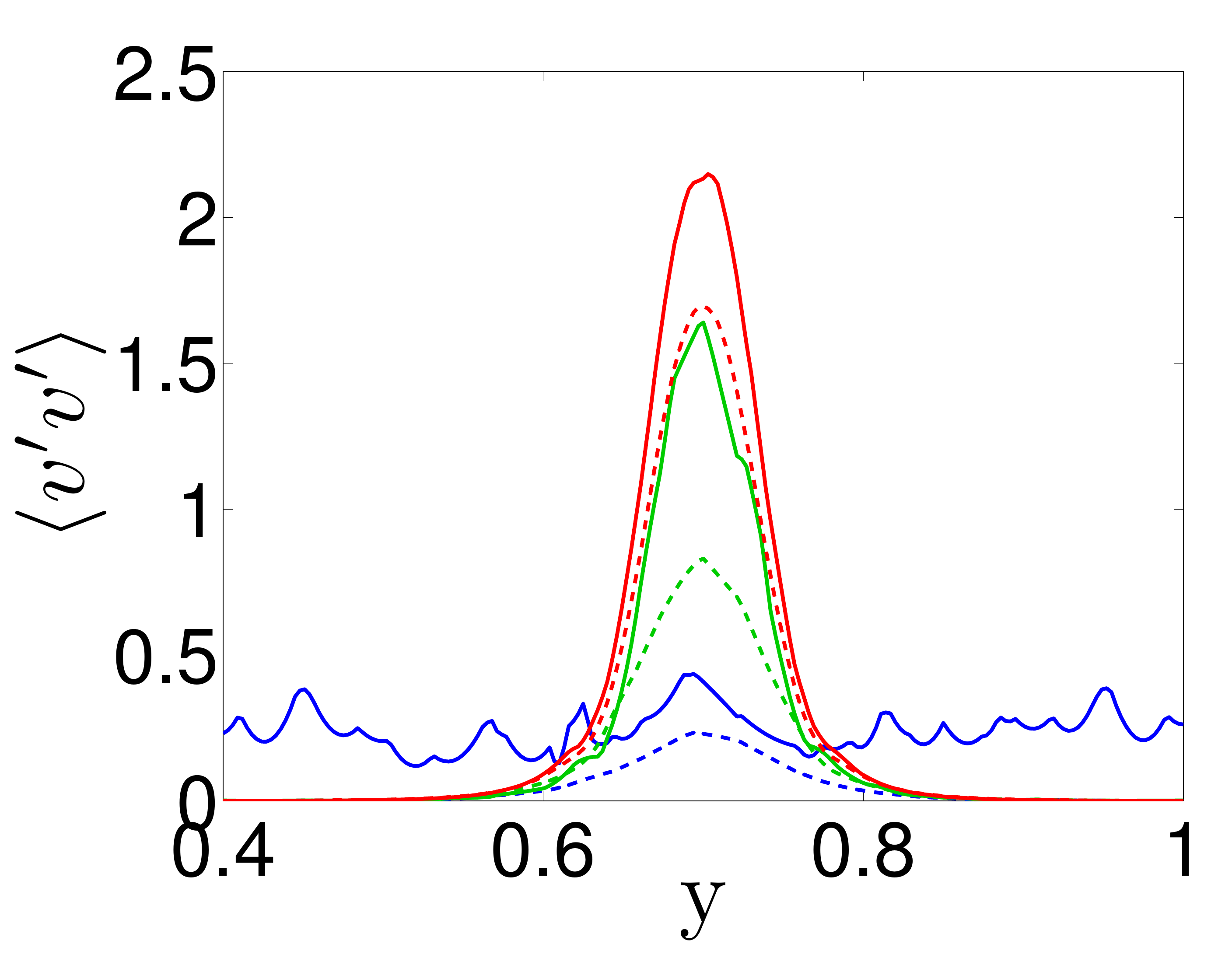}
\includegraphics[height = 1.12in]{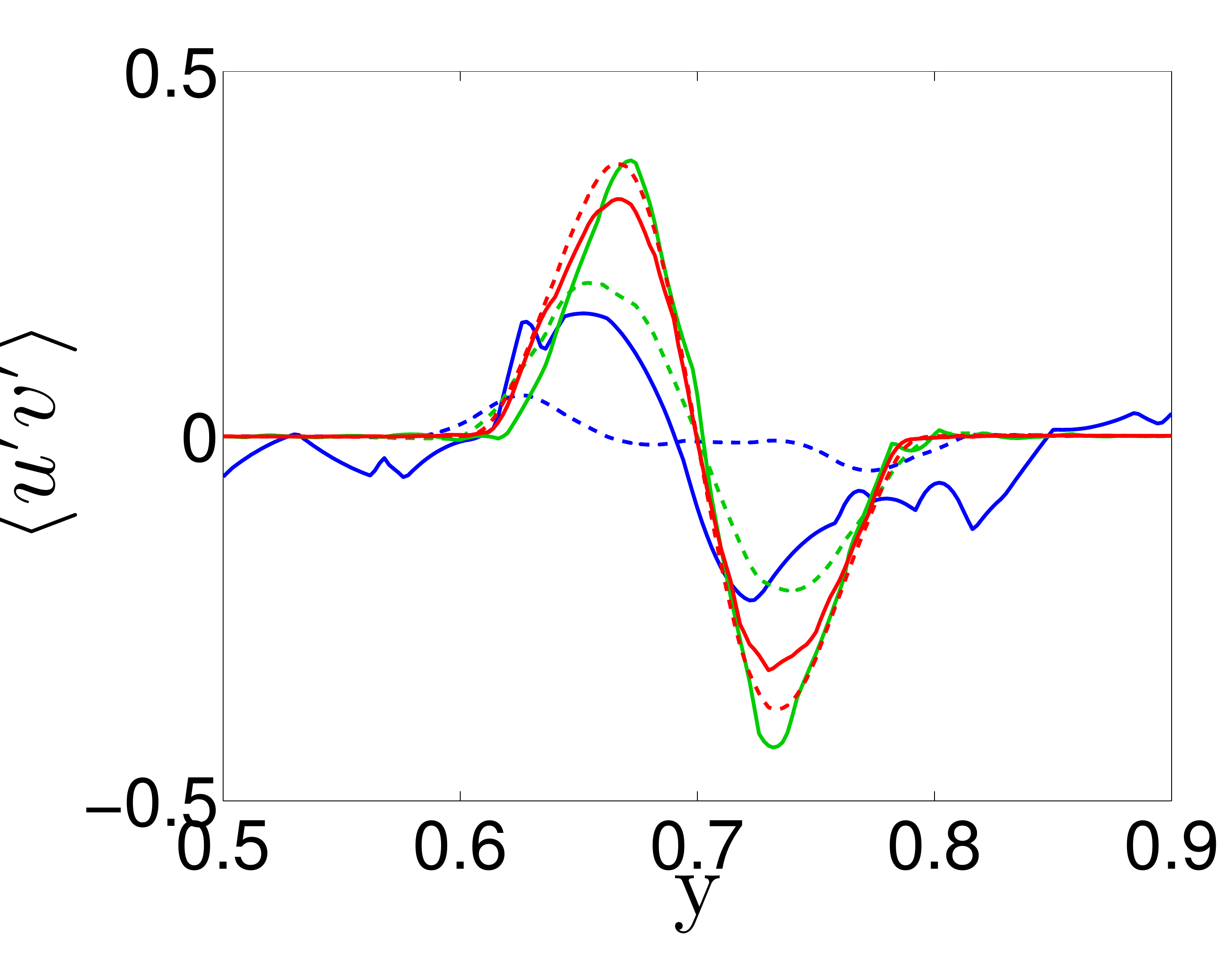}
\\
%\vspace{.1in}
\includegraphics[height = 1.15in]{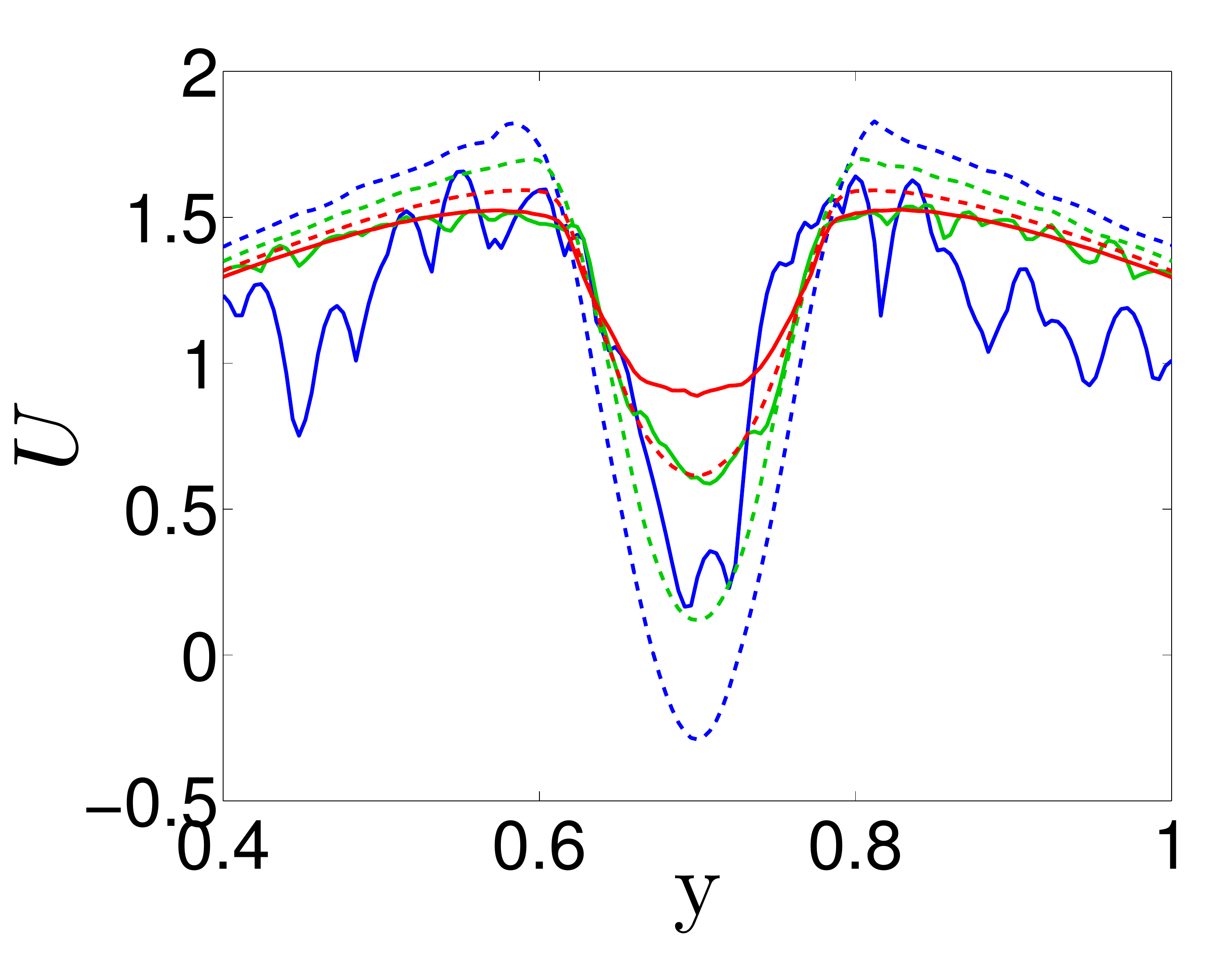}
\includegraphics[height = 1.15in]{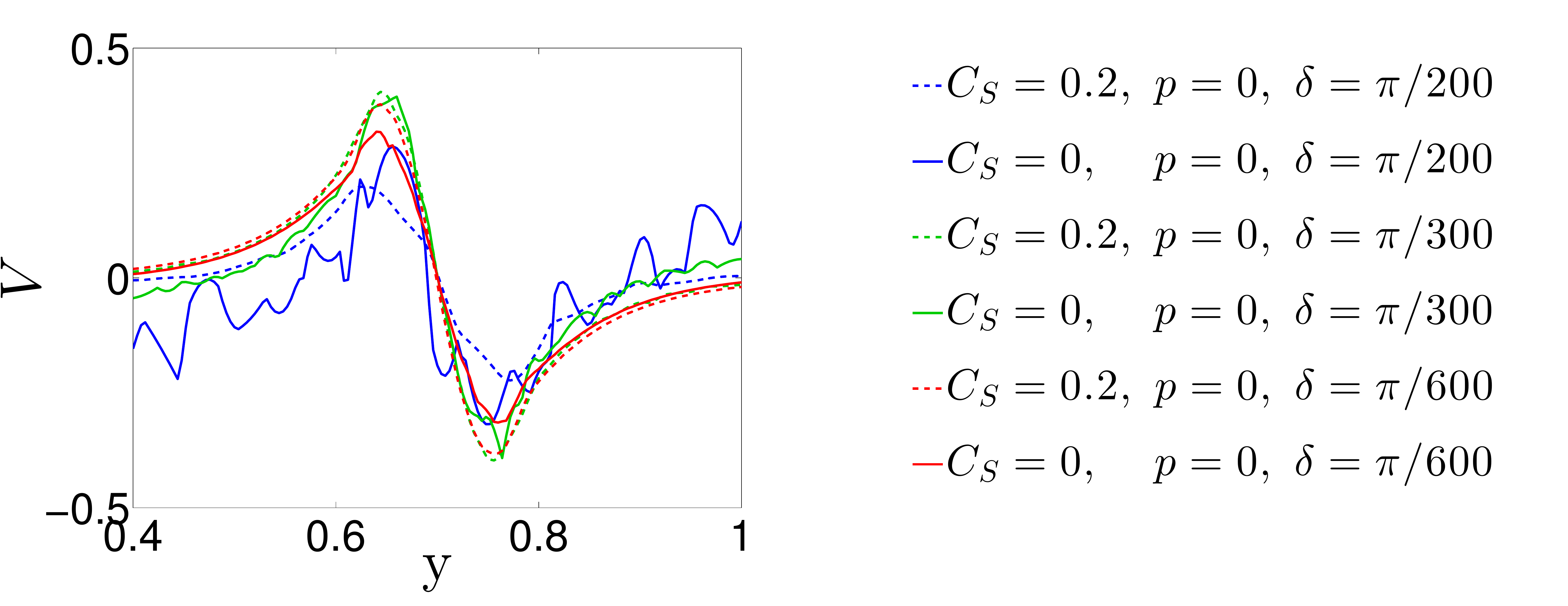}
\caption{Total resolved Reynolds stresses and average velocities along the vertical line at $1D$ downstream for some Smagorinsky models. }
%\vspace{-0.2in}
\label{calib_data}
\end{figure}

The synthetic data are generated by solving Smagorinsky model \eqref{LES-EV}--\eqref{vanDriest} with $C_S=0.15,\, p=0.05$ and $\delta = \pi/480$. The data sets used for calibration process are taken at 11 stations in a distance of $\approx 1D$ downstream. Specifically, these points locate equidistantly on the vertical line $x = 0.65$ between $y = 0.6$ and $y = 0.8 $. For each point, the data of average streamwise and vertical velocities, denoted by $U$ and $V$, as well as total streamwise, vertical and shear Reynolds stresses, i.e., $\langle u'u' \rangle$, $\langle v'v' \rangle$ and $\langle u'v' \rangle$, are selected, giving a total of $55$ reference data. For clarity, the bounds of uniform prior PDFs and the true values of calibration parameters are listed in Table \ref{prior}. In Figure \ref{calib_data}, the measurements of interested velocities and Reynolds stresses along $x=0.65$ are plotted for some typical simulations. We observe that except for $\delta = \pi/200$, the approximated quantities are quite smooth and have expected patterns, see \cite{BM94}. Certainly, the plots show significant differences among different models. In practice, LES models which give distinctly poor results such as those at $\delta = \pi/200$ could be immediately ruled out from the calibration process, informed by the fact that the wall-adjacent grid points lie outside the viscous sublayer. However, it is useful here to examine the response surfaces and the accuracy of the surrogate systems in these cases, and we choose to include these large filter widths in the surrogate domain instead. 

\begin{table}[h] %\vspace{.1in}
\begin{center} %\resizebox{10cm}{!}{
\begin{tabular}{*{4}{c}}
\hline
&  & True value & $\Gamma$ \\ \hline
 Smagorinsky constant & $C_S$ & $0.15$ & $[0,0.2]$ \\ 
 van Driest parameter & $p$ & $0.5$ & $[0,2]$ \\
 Filter width& $\delta$ & $\pi/480$ & $[\pi/600,\pi/200]$ \\
\hline
\end{tabular}
\hspace{0.3cm}\caption{The true parameter values and the initial searching regions for model calibration.}\label{prior}
%\vspace{-.1in}
\end{center}
\end{table}

Finally, it is worth mentioning that Smagorinsky model coefficients are not the only parameters that influence the quality of LES solutions. Indeed, other numerical parameters such as time step size and averaging time also have significant impacts, see, e.g., \cite{RFBP97,Bre98}. While an estimation of their influence is not conducted here, we need to ensure that the errors caused by them do not dominate the uncertainties in the calibration parameters. A simple validation test is carried out on the flow statistics generated by Smagorinsky model of $C_S=0,\, p = 0$ and $\delta = \pi/480$. The flow simulation is replicated first with the temporal resolution refined by a factor of two, i.e., $\Delta t = 0.005$, and then with a doubled averaging period, i.e., by setting $T = 19$. We also conduct another simulation in which the zero gradient replaces Dirichlet outflow boundary condition to justify that the numerical oscillation at the downstream boundary does not disturb the inner domain. The maximum change in five velocity and Reynold stress profiles of interest in these modified models is presented in Table \ref{valid}. We see that among three investigated source of numerical errors, the temporal resolution is the most prominent, as it makes up approximately $80\%$ of the change in all data. More importantly, Table \ref{valid} reveals that the total maximum change is approximately $0.05$ in the vertical Reynolds stress data and $0.025$ for other quantities. Numerical errors of the synthetic calibration data, as well as model outputs, are expected to be around these values. In the uncertainty analysis following, for the \eqref{MVN} likelihood model, we will assume that the reference data are corrupted with Gaussian random noise of $0.1$. 

\begin{table}[h]

\centering 
%\resizebox{10cm}{!}{
\begin{tabular}{*{6}{c}}
\hline
Component modified & $U$ & $V$ & $\langle u'u' \rangle$ & $\langle v'v'\rangle$ & $\langle u'v' \rangle$  \\ \hline
Time step  & $0.0148$ & $0.0213$ & $0.0202$ & $0.0371$ & $0.0186$ \\
Averaging period & $0.0048$ & $0.0025$ & $0.0014$ & $0.0011$ & $0.0045$ \\
Outflow BC & $0.0022$ & $0.0029$ & $0.0016$ & $0.0061$ & $0.0020$ \\
\hline
\end{tabular}
\caption{The maximum change in average velocity and Reynold stress profiles under the modifications of time step, averaging period and outflow BC. }
%\vspace{-.1in}
\label{valid}
\end{table}

\subsection{Results and discussions}
\label{sec:result}
This section justifies the accuracy and efficiency of the surrogate modeling method described in \S\ref{sec:surrogate}, when applied to the numerical example of two-dimensional flow around a cylinder specified in \S\ref{sec:les}. We utilize the software package \textit{FreeFem++} \cite{Hec12} in solving the Smagorinsky discretization scheme. The adaptive sparse grid interpolation and integration schemes are generated using functions in the TASMANIAN toolkit \cite{Sto13}. The DRAM algorithm \cite{HLMS06} is chosen for MCMC sampling of the surrogate PPDF. 

The surrogate system for outputs is constructed using the linear basis functions, first on the standard sparse grid of level 5, then the grids are refined adaptively up to level 8. The total numbers of model executions needed for the four interpolants are $177,\, 439,\, 1002$ and $2190$, respectively, which are also the number of points of the four corresponding adaptive sparse grids. 

\begin{figure}[h]
\centering
\includegraphics[height = 1.7in]{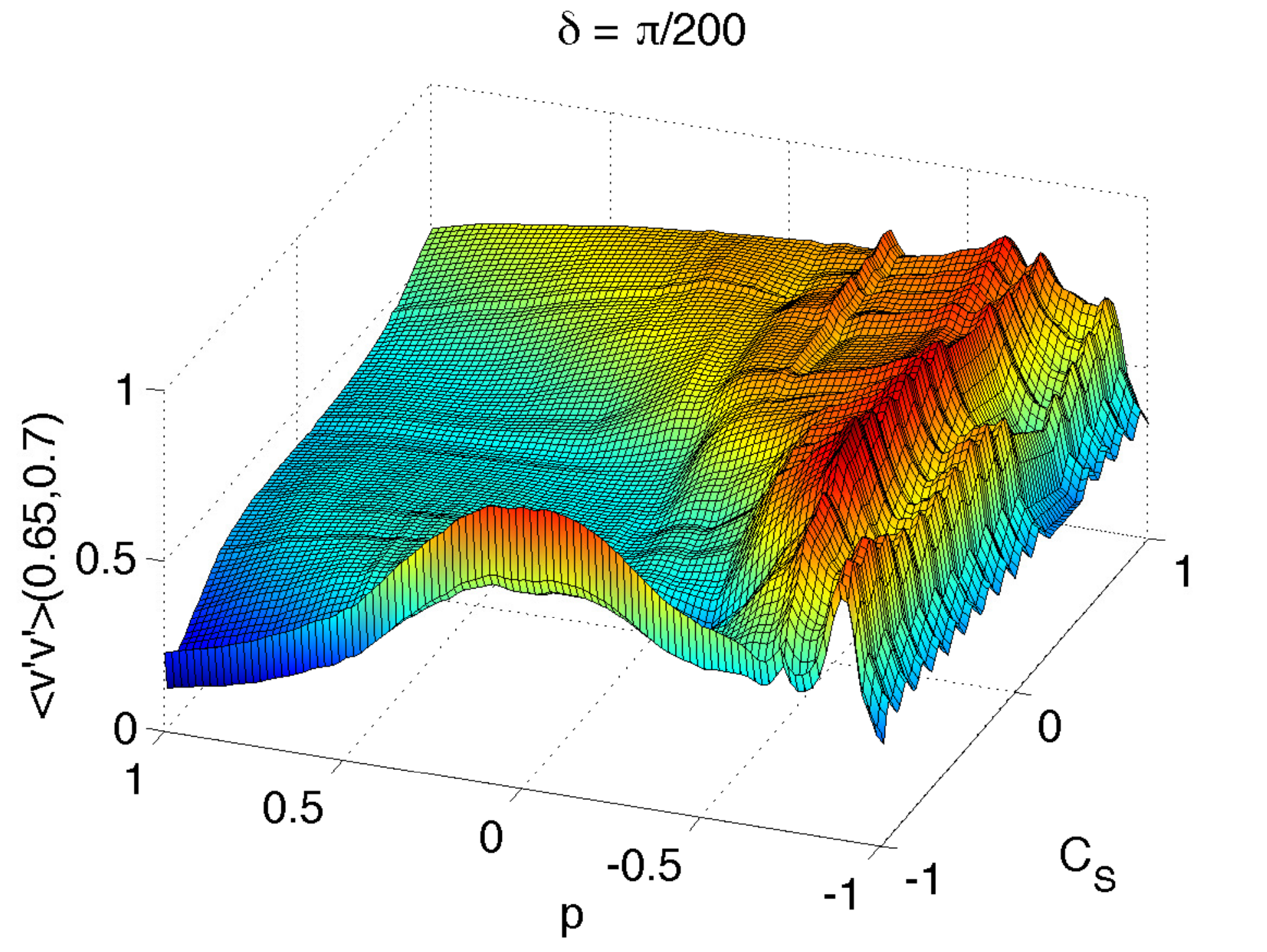}
\includegraphics[height = 1.7in]{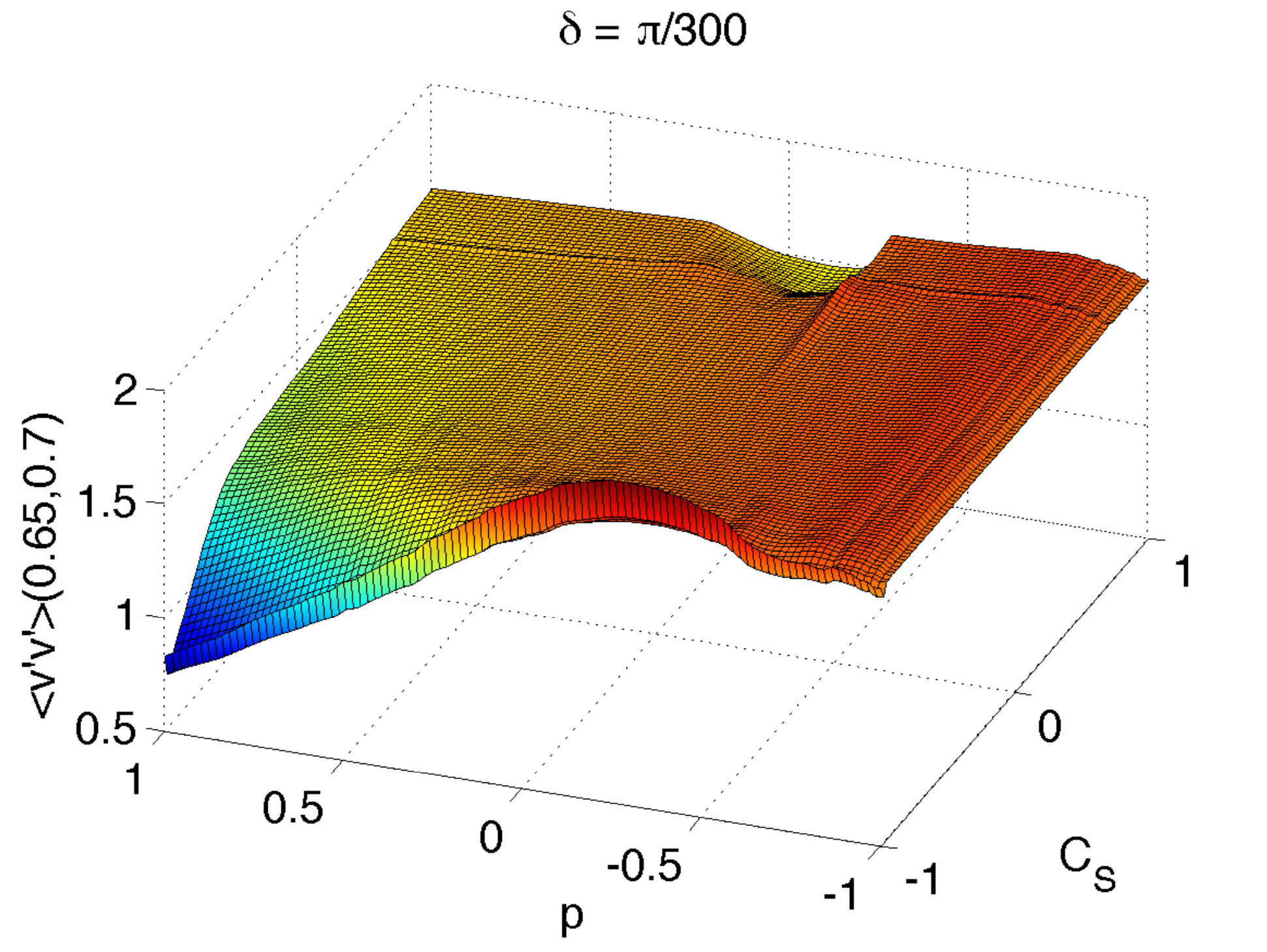}
\includegraphics[height = 1.7in]{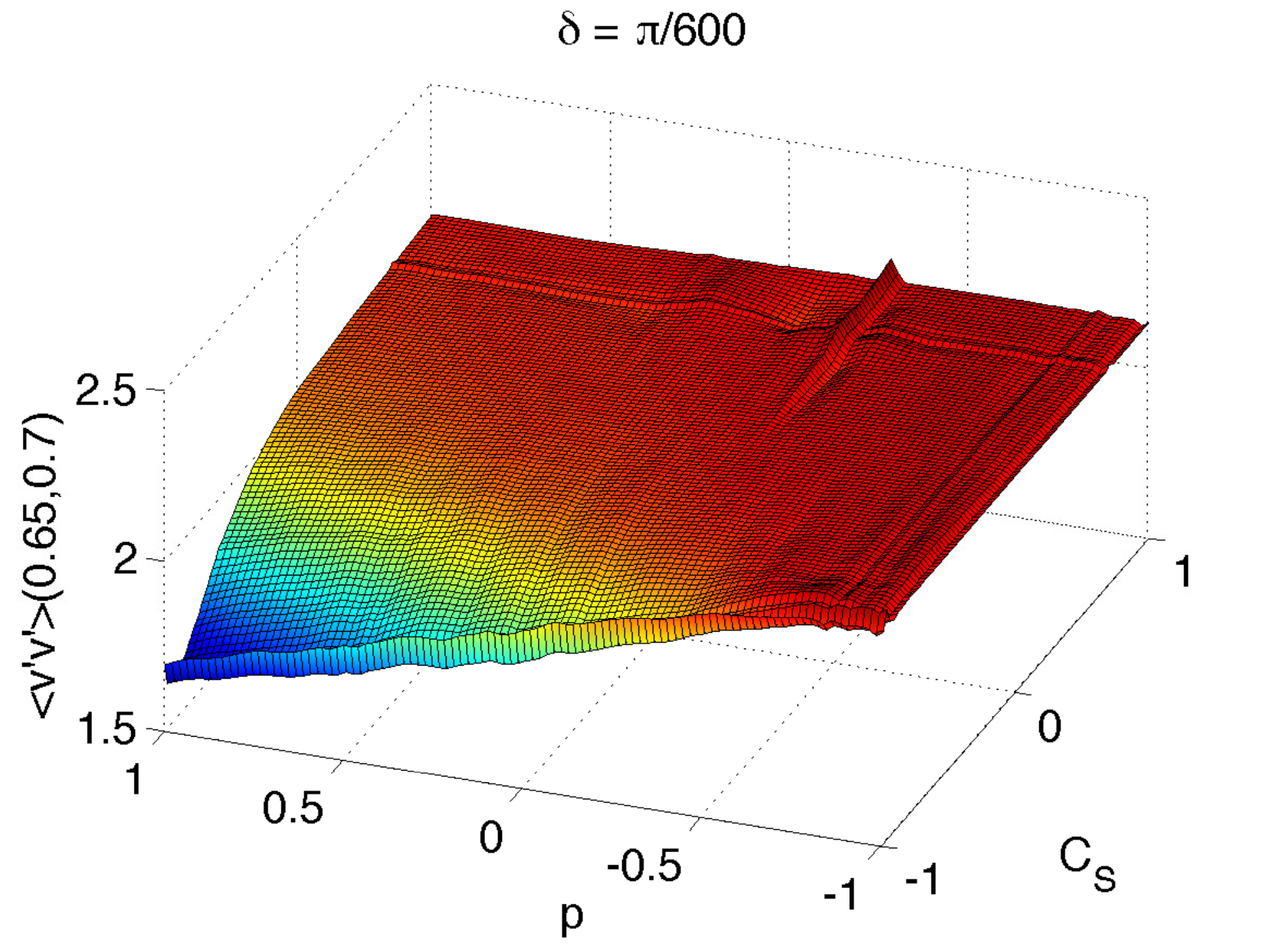}
\caption{Surfaces of the predicted vertical Reynolds stress data at $ (0.65,0.7)$ generated by the AHSG method at level $8$. $C_S$ and $p$ are normalized such that their searching regions are $[-1,1]$.}
\label{surface}
\end{figure}

\begin{figure}[h]
\centering
\includegraphics[height = 1.7in]{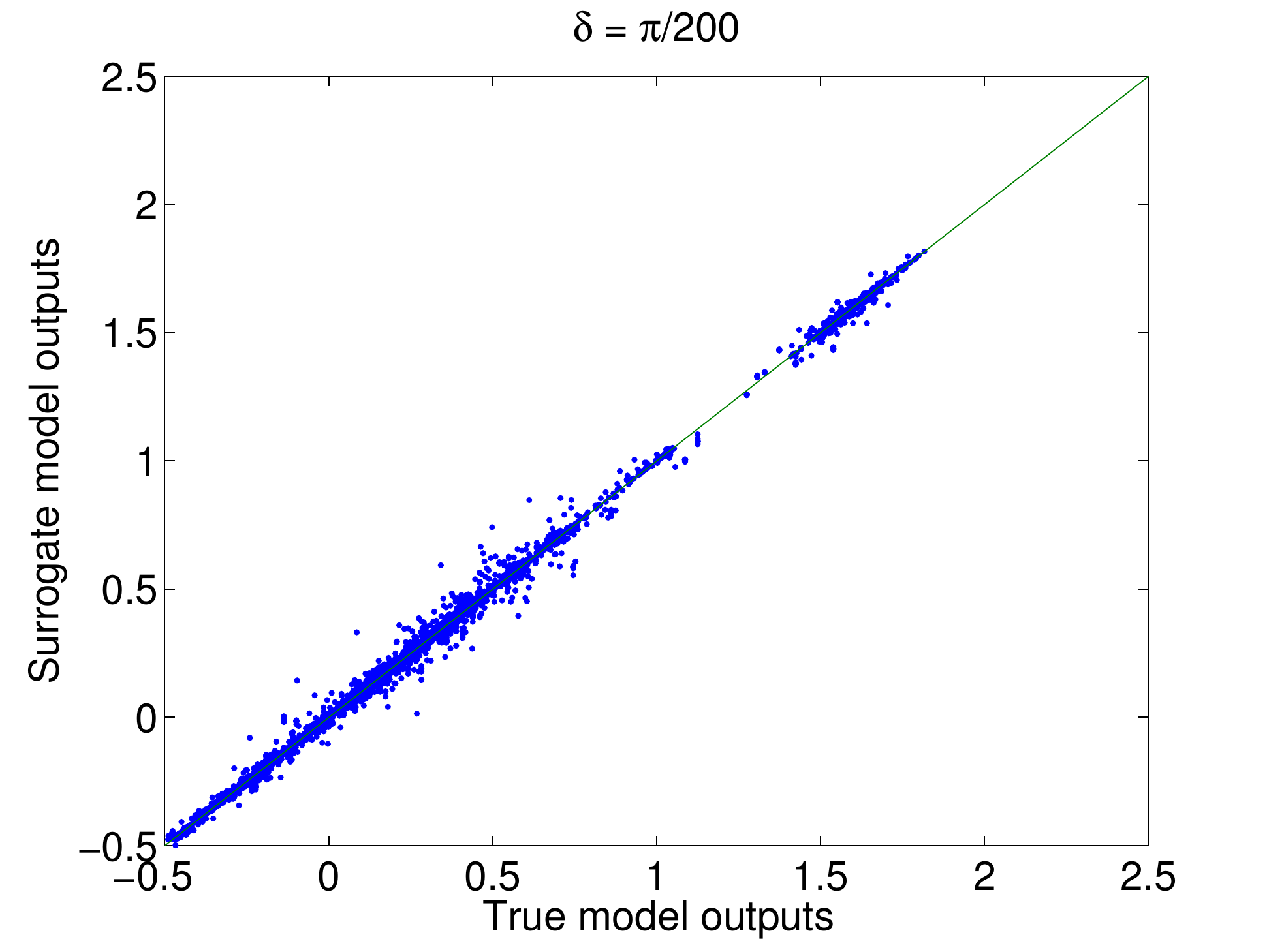}
\includegraphics[height = 1.7in]{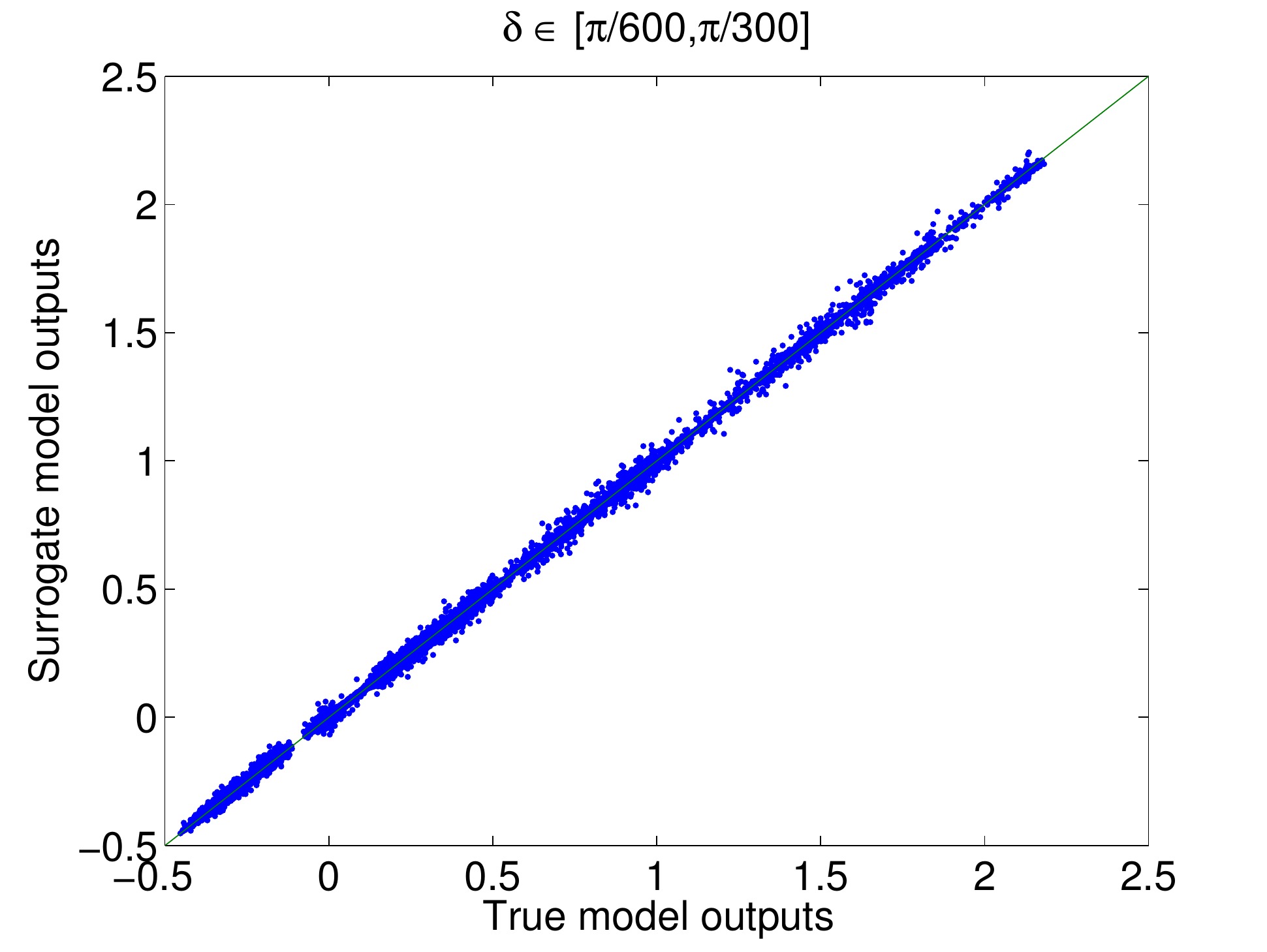}
\caption{Scatter plots for the prediction of the output data given by the surrogate system on level $7$ sparse grid.}
\label{true_vs_surrogate}
\vspace{-.1in}
\end{figure}

The accuracy of a surrogate modeling approach based on the AHSG method is largely determined by the smoothness of the surrogate system, so it is worth examining the surface of the output data in the parameter space. For brevity, we only plot here the vertical Reynolds stress data at the centerline, i.e., $\langle v'v' \rangle (0.65,0.7)$, which is among the most fluctuating (See Figure \ref{calib_data}). Figure \ref{surface} represents some surfaces for typical values of filter width generated on level 8 grid. We observe that the surface according to $\delta = \pi/200$ differs from two other cases ($\delta = \pi/600,\, \delta = \pi/300$) that are remarkably rougher. This, together with Figure \ref{calib_data}, confirms the connection between the complexity of the output function in both the 
physical and parameter spaces. In Figure \ref{true_vs_surrogate}, the scatter plots for the predicted outputs obtained with the surrogate system at level $7$ are presented. The approximations show clear improvement in accuracy with $\delta\in [\pi/600,\pi/300]$, compared to those at larger values. While not considered herein, it is reasonable to expect that the surrogate outputs at least maintain the same accuracy for $\delta \le \pi/600$, since more grid refinement will remove extra non-physical wiggles. In the next part, we justify that this level of accuracy is sufficient for our surrogate-based MCMC method. Although the surrogate systems show remarkable discrepancy for large $\delta$, as previously mentioned, these values, leading to visibly inadequate outputs, should be excluded in practical calibration processes. While the original domain of $\delta$ is $[\pi/600,\pi/200]$, by choosing its true value as $\pi/480$, the effective searching region of $\delta$ is restricted to $[\pi/600,\pi/300]$. 

%Figure \ref{PPDF1}--\ref{PPDF2} shows the marginal posterior PDFs for $C_S$, $m$ and $\alpha$. The components of random vector $\varepsilon$ are taken to be independent, Gaussian random variables with mean zero and prescribed standard deviations. In particular, the standard deviations are set to 2\% of the reference value for the velocity measurements and 5\% of the reference value for the Reynolds stress measurements. 

%\begin{figure}[!h]
%\centering
%\includegraphics[height=5.5cm]{TKElateral_1-eps-converted-to.pdf} 
%\includegraphics[height=5.5cm]{TKEcentre-eps-converted-to.pdf} 
%\caption{Total resolved turbulent kinetic energy at $x/D = 1$ (left) and along the centerline (right) for different EV models and filter radii. }
%\label{energy}
%\end{figure}

\begin{figure}[h]
\centering
\includegraphics[height=3.cm]{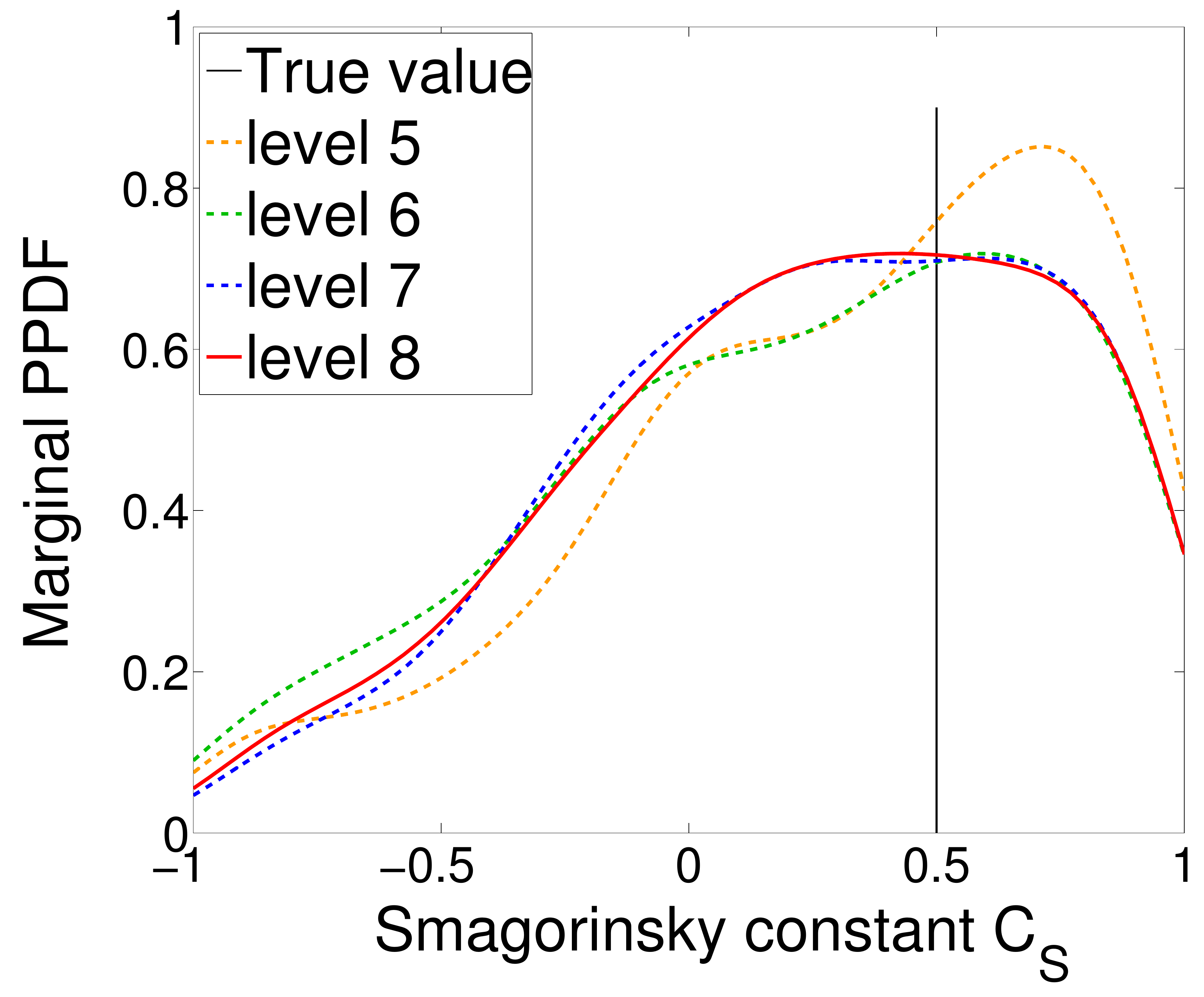} 
\includegraphics[height=3.cm]{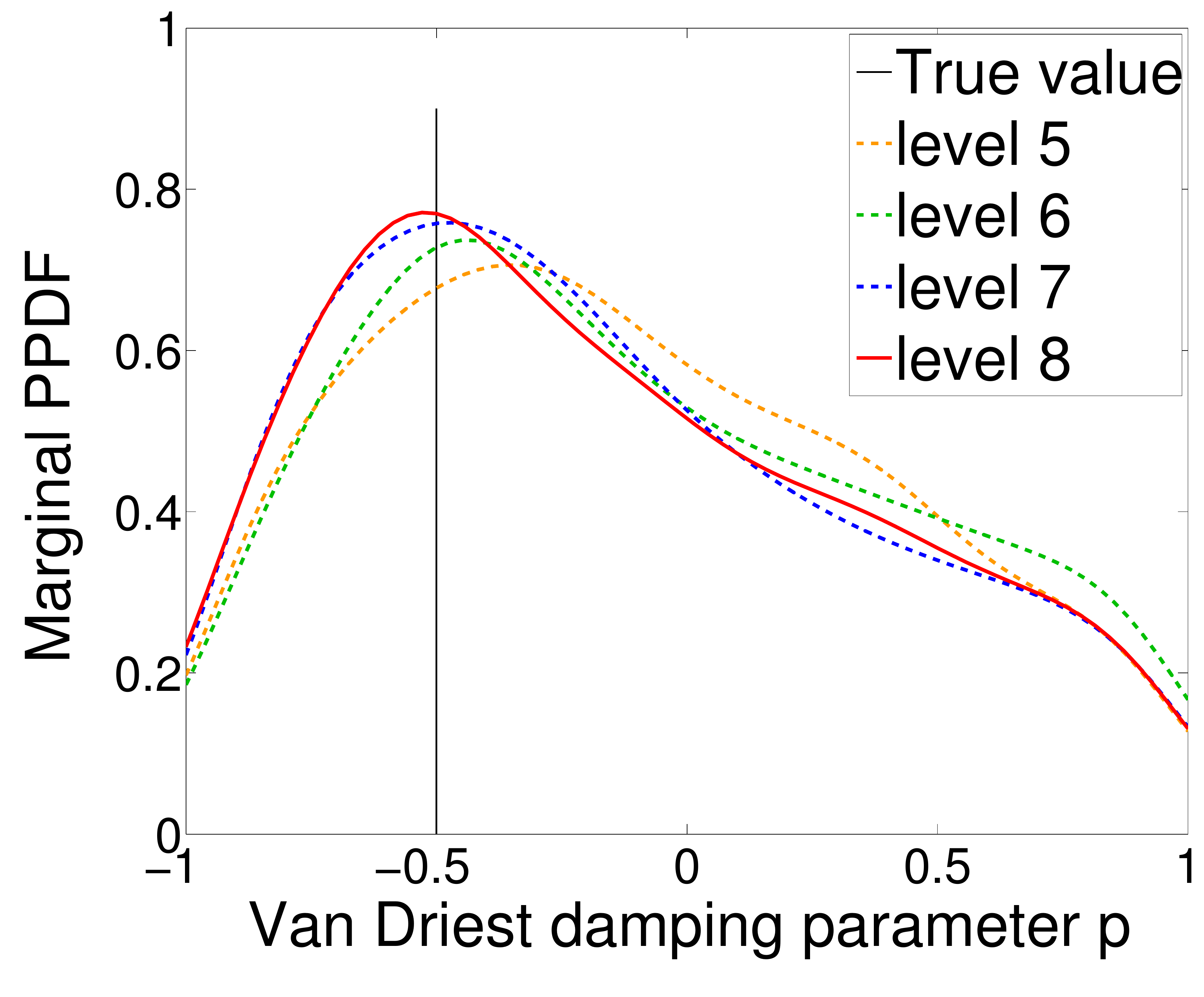} 
\includegraphics[height=3.cm]{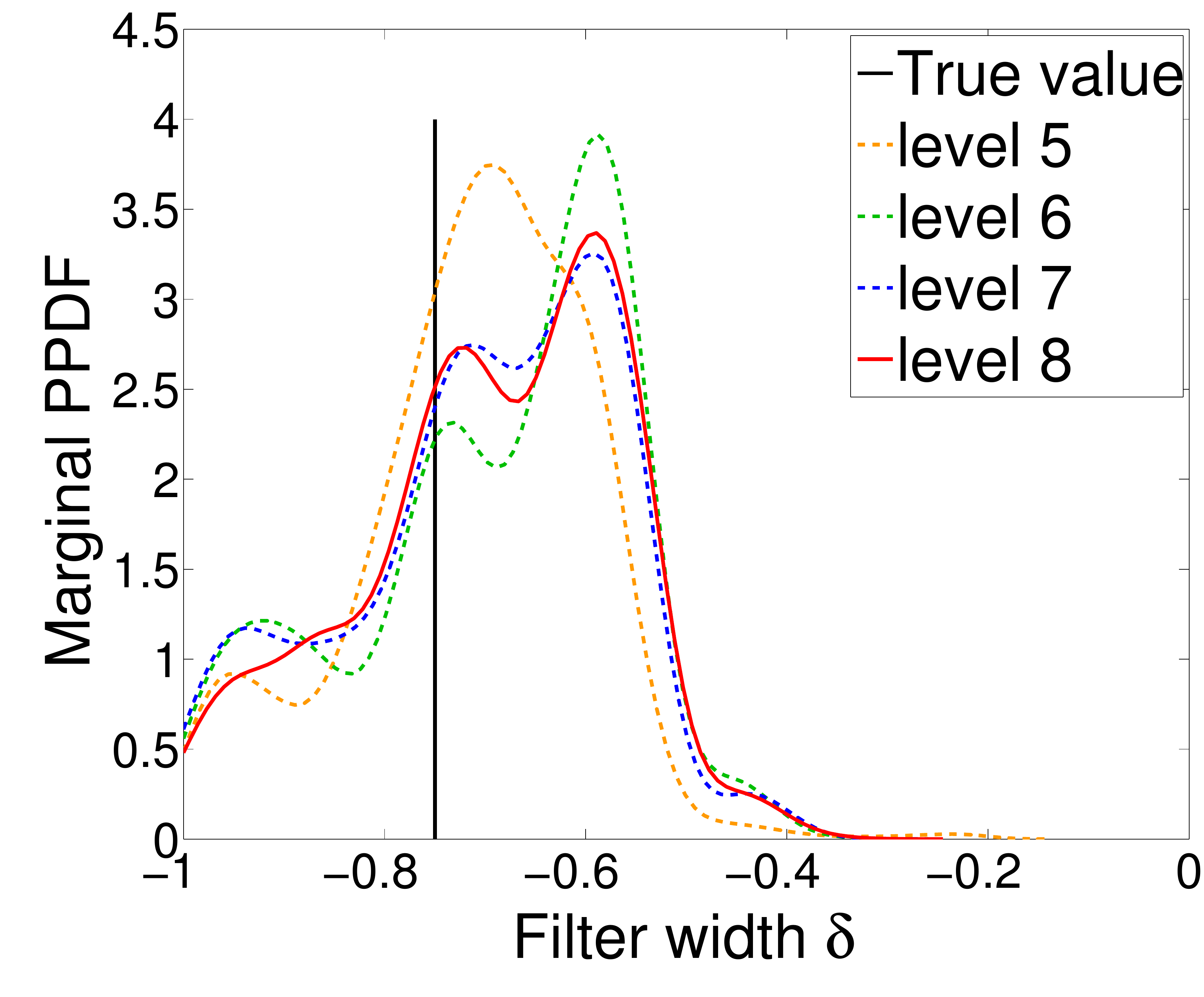} 
\caption{Marginal posterior probability density functions of three Smagorinsky model parameters with \eqref{MVN} likelihood model estimated using the linear surrogate systems on level $5$--$8$ adaptive sparse grids.}
\label{fig:PPDF1}
\vspace{-.1in}
\end{figure}

To evaluate the accuracy and efficiency of our surrogate modeling approach, the DRAM-based MCMC simulations using the surrogate PPDF $\tilde{P}(\bm{\theta}|\bm{d})$ constructed in Algorithm \ref{surrogate_construct} are conducted. Each MCMC simulation draws 60,000 parameter samples, the first 10,000 of which are discarded and the remaining 50,000 samples are used for estimating the PPDF.  For the first experiment, \eqref{MVN} likelihood function is employed; the data are corrupted by $10\%$ Gaussian random noise, treated as numerical errors. Figure \ref{fig:PPDF1} plots the marginal PPDFs where the three parameters are normalized such that the searching region is $[-1,1]^3$. The black vertical lines represent the true values listed in Table \ref{prior}. The red solid lines are the marginal PPDFs estimated by MCMC simulations based on the surrogate systems on level 8 grid, and the dashed lines represent those based on the surrogate systems on lower levels. The figure indicates that the MCMC results according to level 7 and level 8 sparse grids, which require 1002 and 2190 model executions correspondingly, are already close to each other. Thus, the surrogate PPDF on level 8 is accurate enough for MCMC simulations.

\begin{figure}[h]
\centering
\includegraphics[height=3.cm]{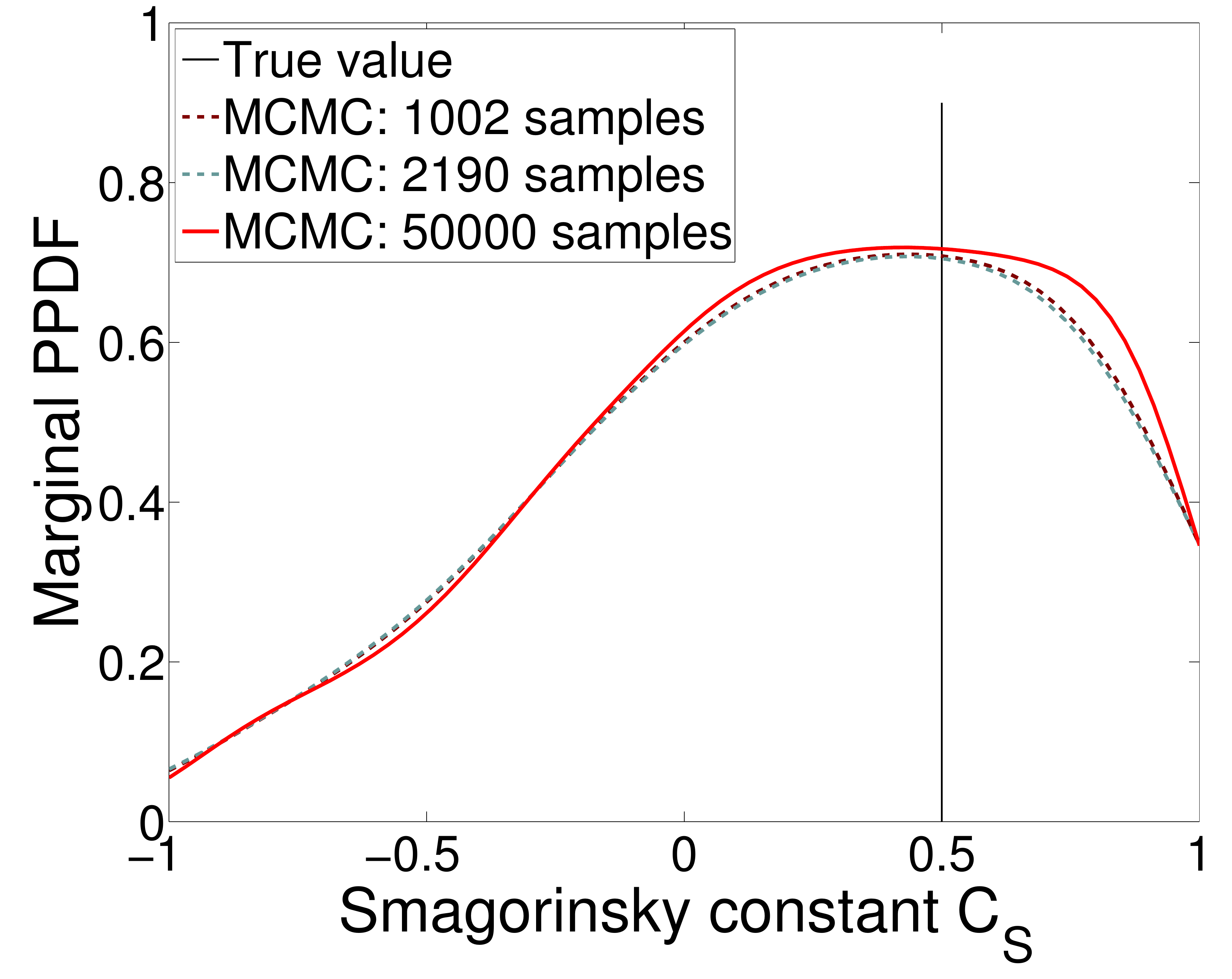} 
\includegraphics[height=3.cm]{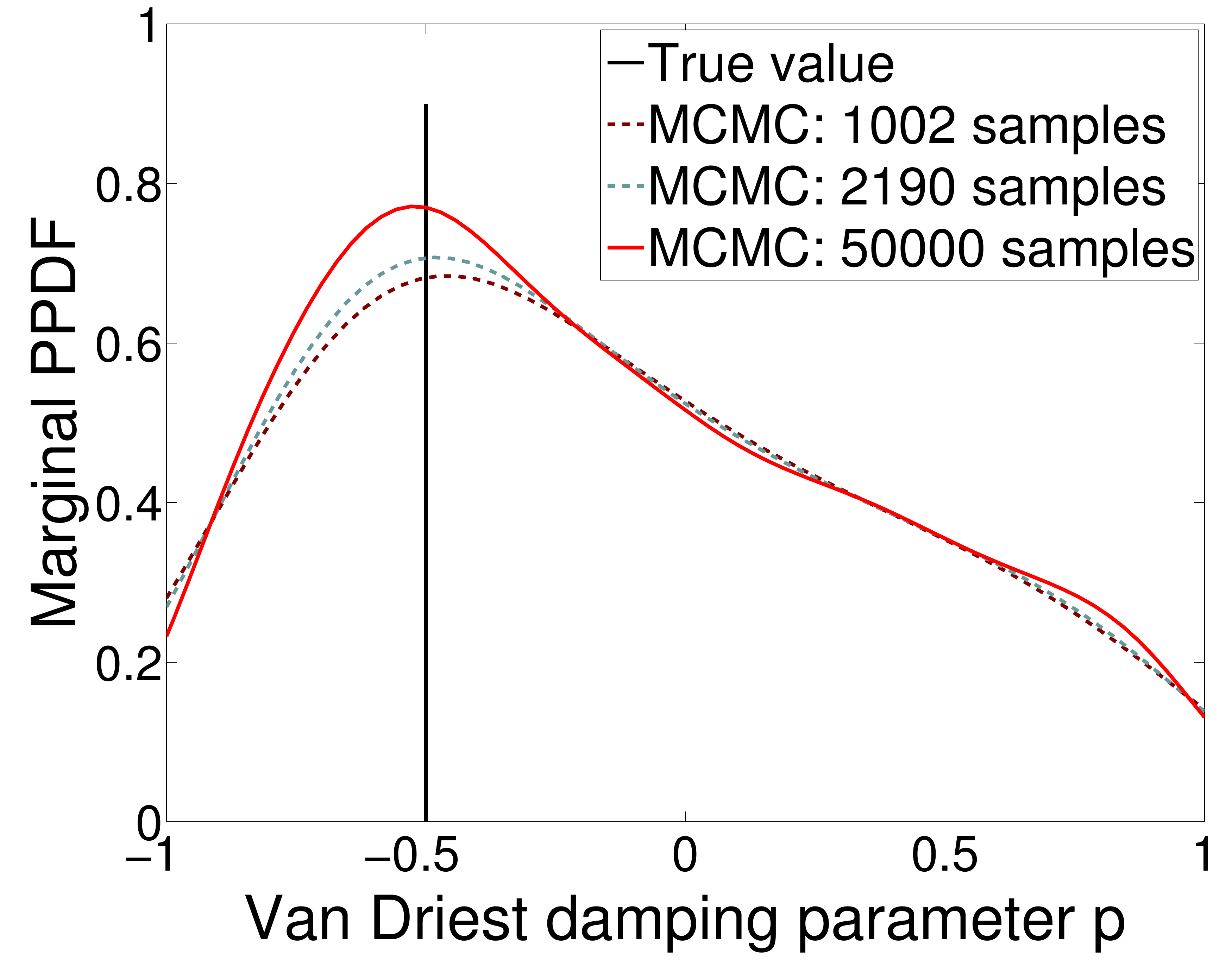} 
\includegraphics[height=3.cm]{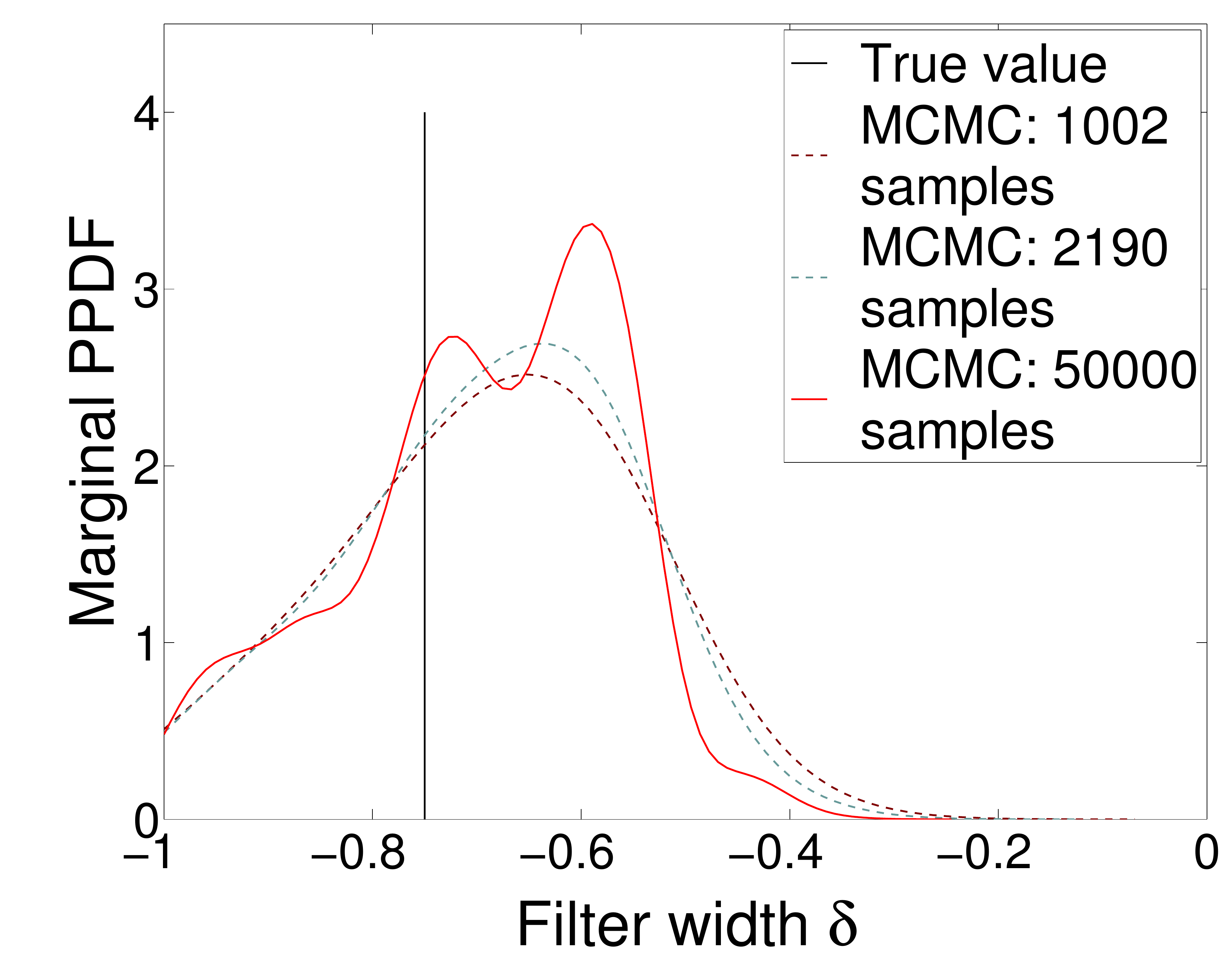} 
\caption{Marginal posterior probability density functions of model parameters with \eqref{MVN} likelihood function estimated using the linear surrogate systems on level 8 adaptive sparse grids with 1002, 2190 and 50000 samples {(excluding 10000 samples for burn-in period). These are the numbers of {model executions} that the conventional MCMC requires to obtain similar results.}}
\label{fig:PPDF2}
\vspace{-.1in}
\end{figure}

{We proceed to compare the accuracy of the surrogate-based with the conventional MCMC with equal computational effort, i.e., same number of model executions. Due to the high computational cost, a proper conventional MCMC simulation is not conducted in this work. However, given the accuracy of the surrogate system, we expect that marginal PPDFs obtained from conventional MCMC are very close with those from surrogate-based MCMC on high-level grid and therefore, run the MCMC simulation with samples drawn from level 8 surrogate. The first 10000 samples are discarded to minimize the effect of initial values on the posterior inference. Figure \ref{fig:PPDF2} depicts the marginal PPDFs for model parameters obtained with 1002, 2190 and 50000 samples after burn-in period. Let us remark that if conventional MCMC is employed, these are the numbers of \textit{model executions} required to obtain similar results.} %Should the conventional MCMC be applied, the same numbers of model executions would be required for producing similar PPDFs, not to mention executions needed for initialization of the Markov chain. 
Comparing Figure \ref{fig:PPDF1} and \ref{fig:PPDF2} indicates that with the same number of model executions, the approximations using surrogate system are more accurate than those using conventional MCMC, highlighting the efficiency of our surrogate modeling method. 

\begin{figure}[h]
\centering
\includegraphics[height=3.cm]{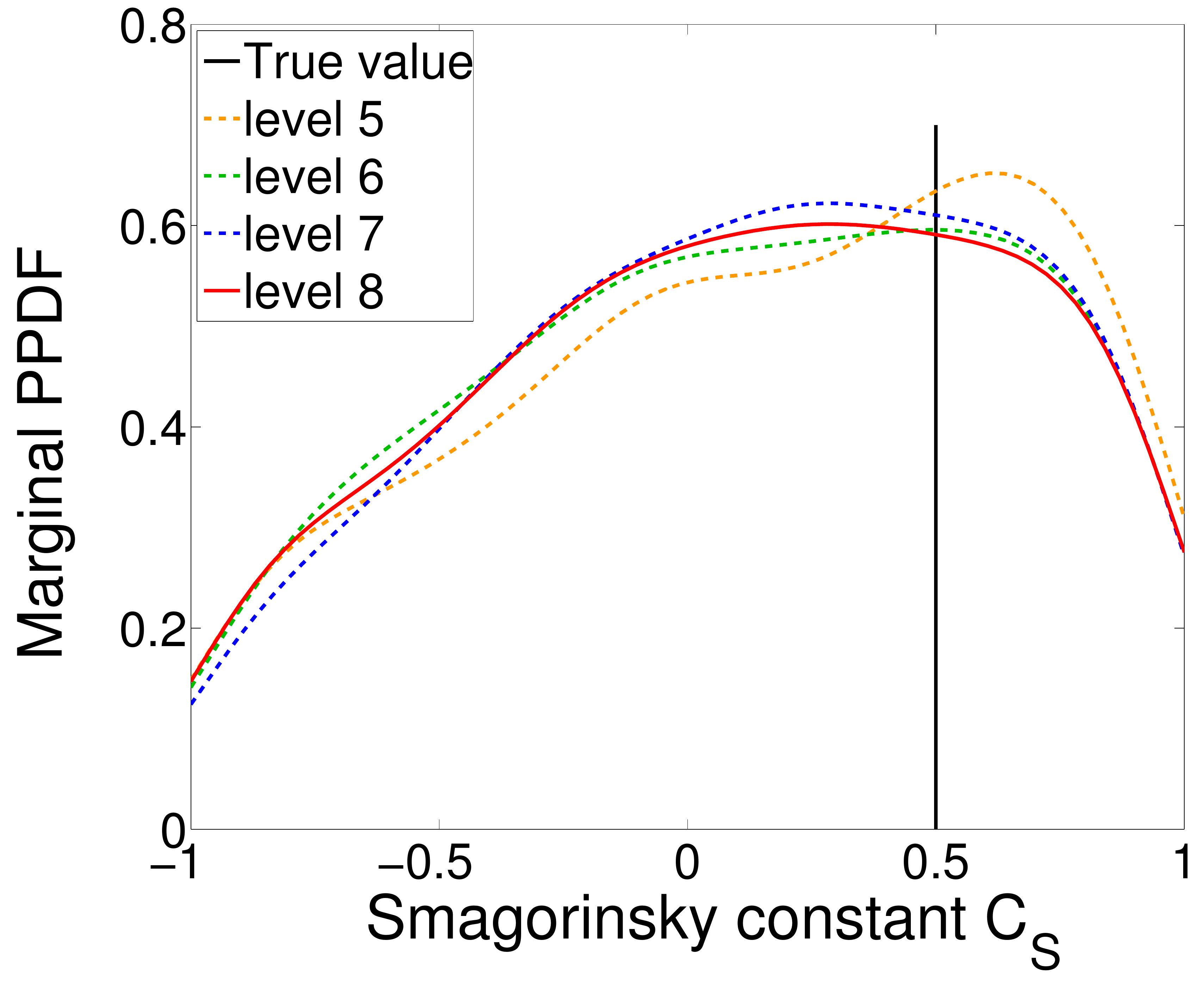} 
\includegraphics[height=3.cm]{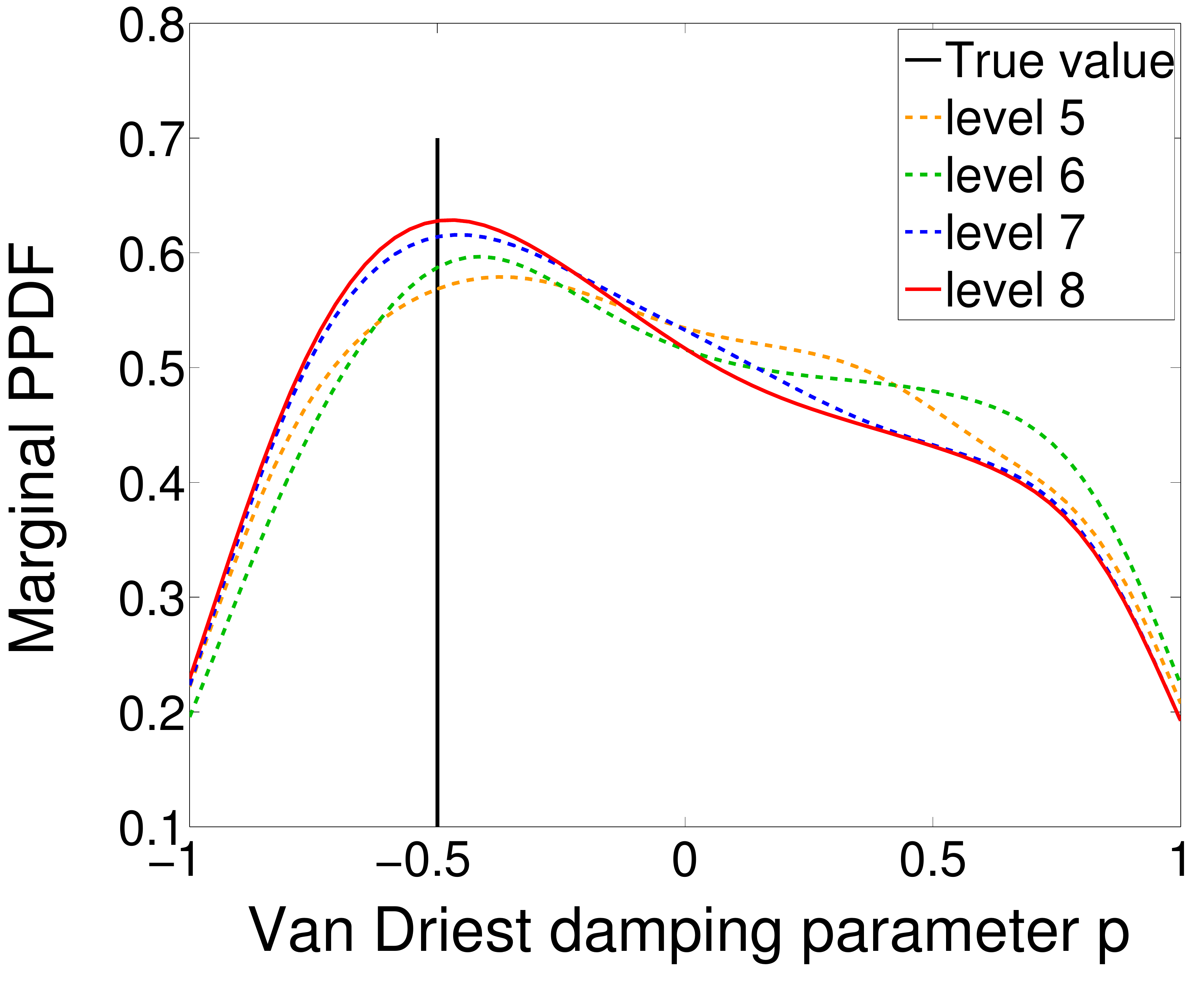} 
\includegraphics[height=3.cm]{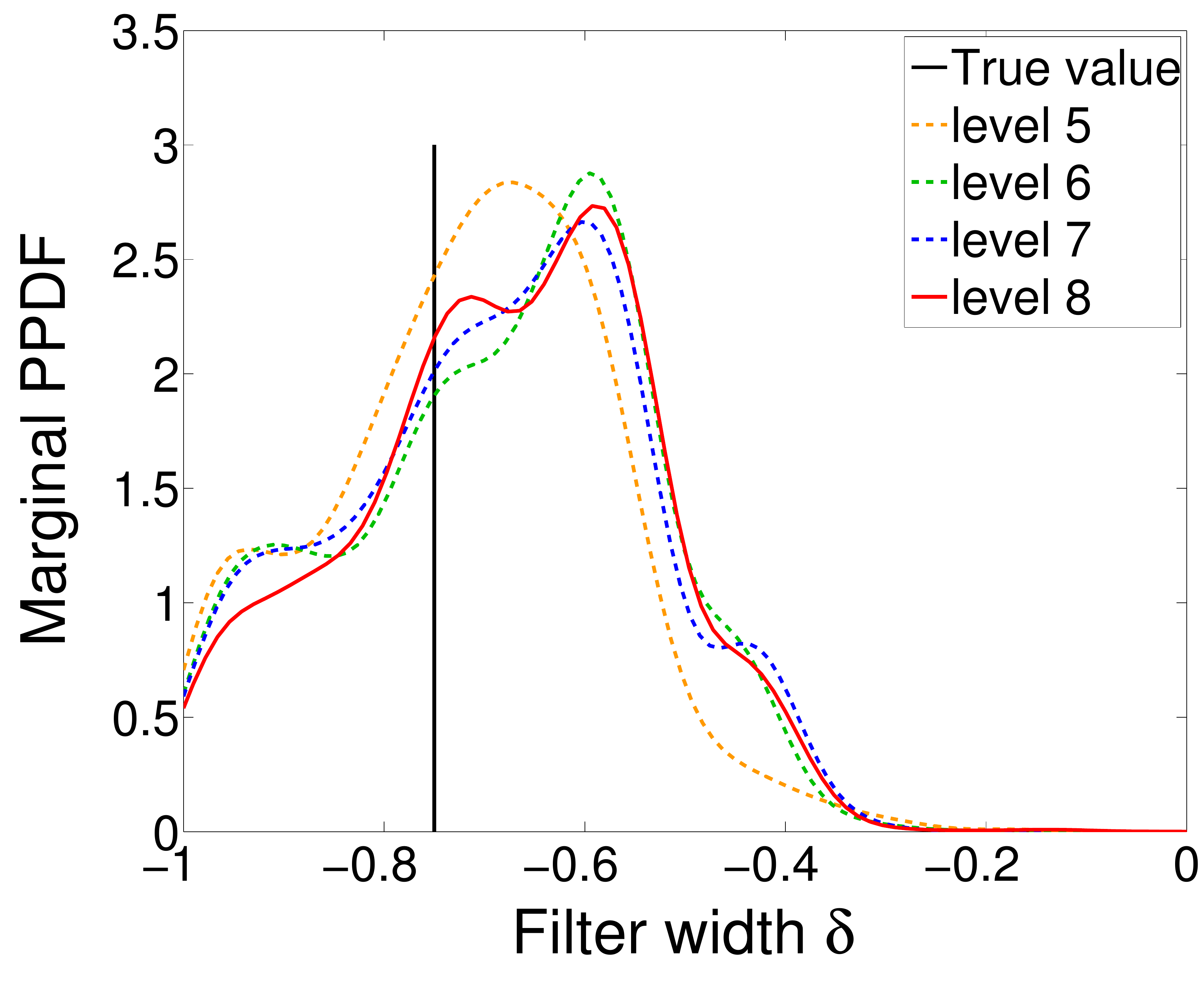} 
\caption{Marginal posterior probability density functions of three Smagorinsky model parameters with \eqref{EXP} likelihood model estimated using the linear surrogate systems on level $5$--$8$ adaptive sparse grids.}
\label{fig:PPDF3}
%\vspace{-.1in}
\end{figure}

\begin{figure}[h]
\centering
\includegraphics[height=3.cm]{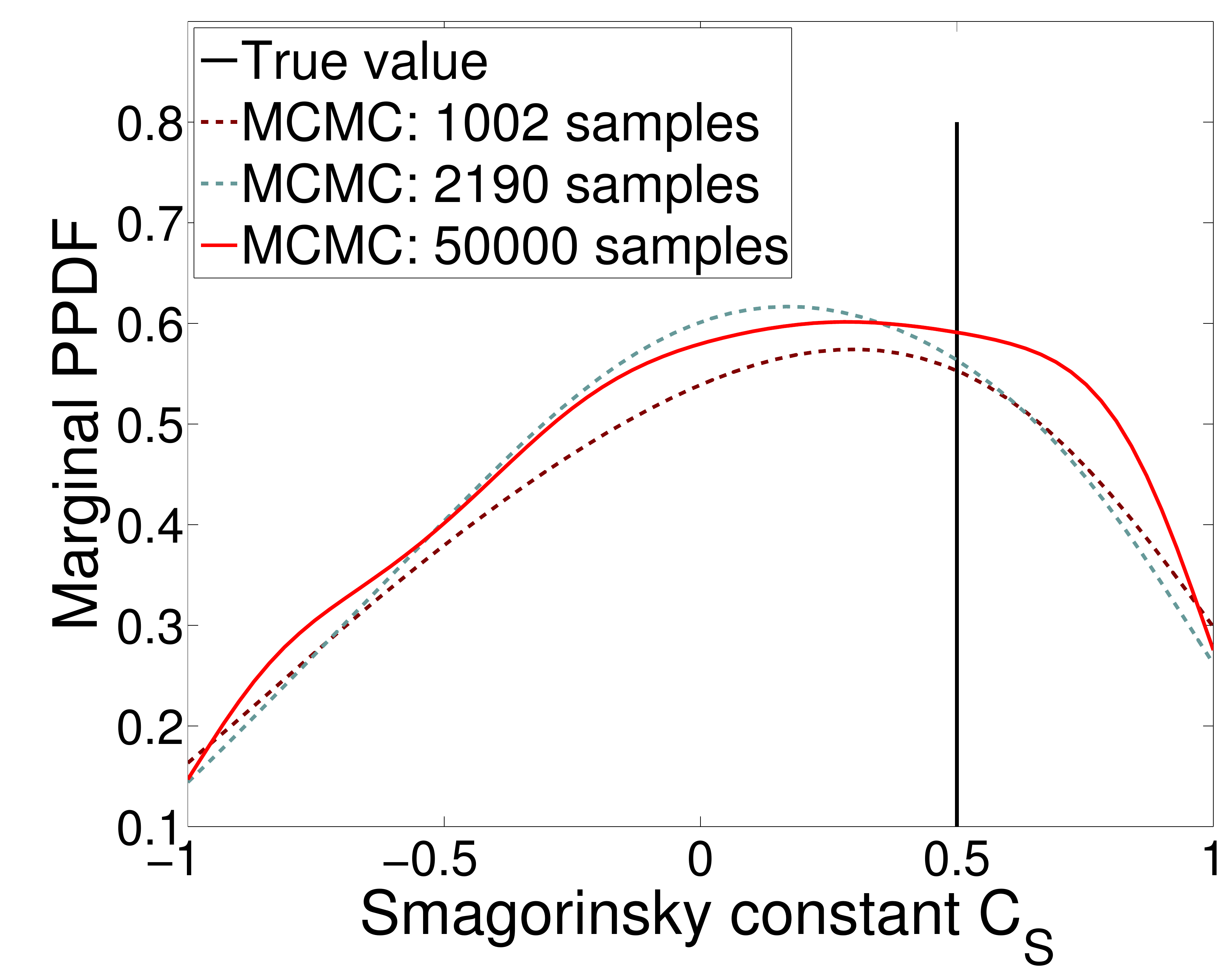} 
\includegraphics[height=3.cm]{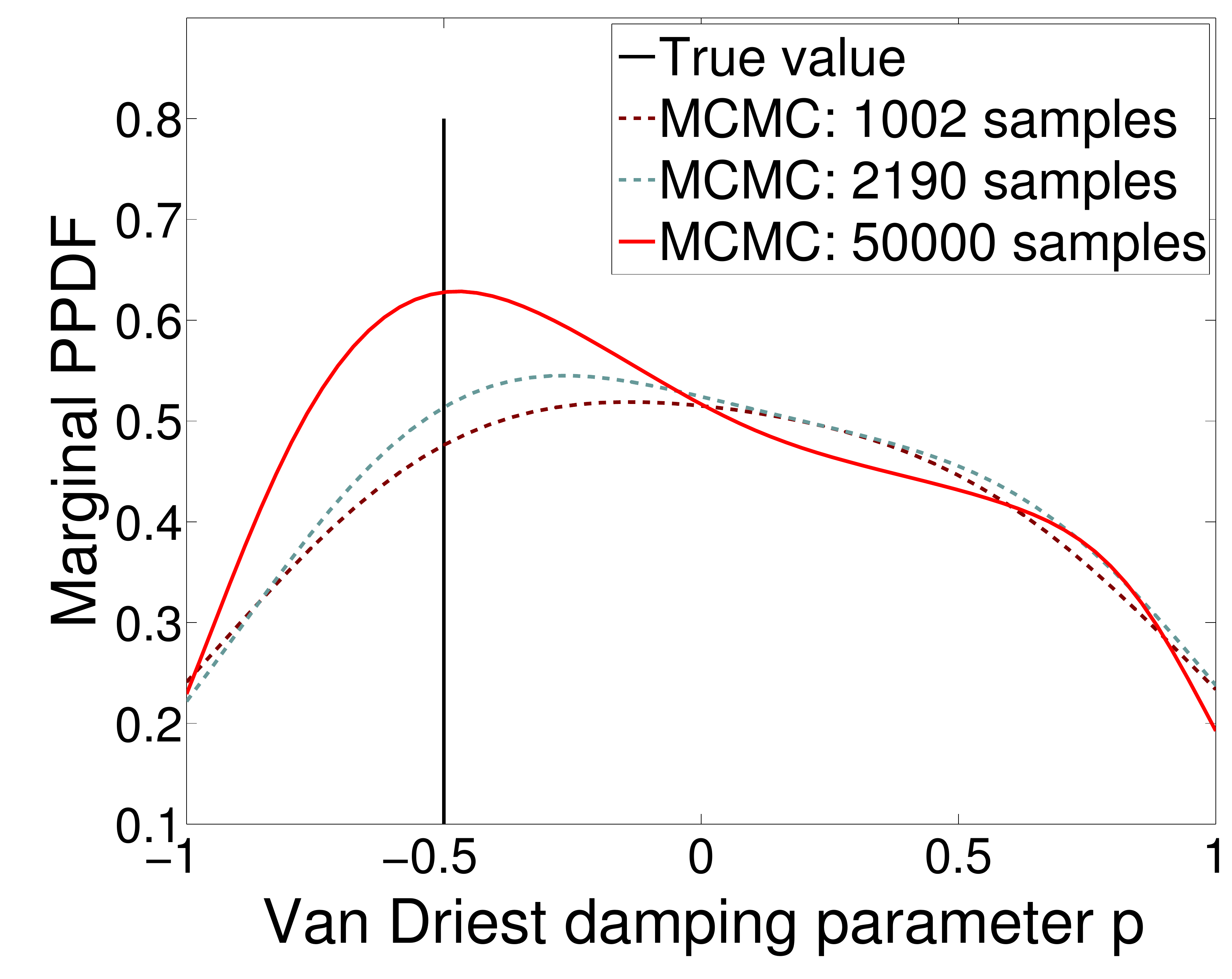} 
\includegraphics[height=3.cm]{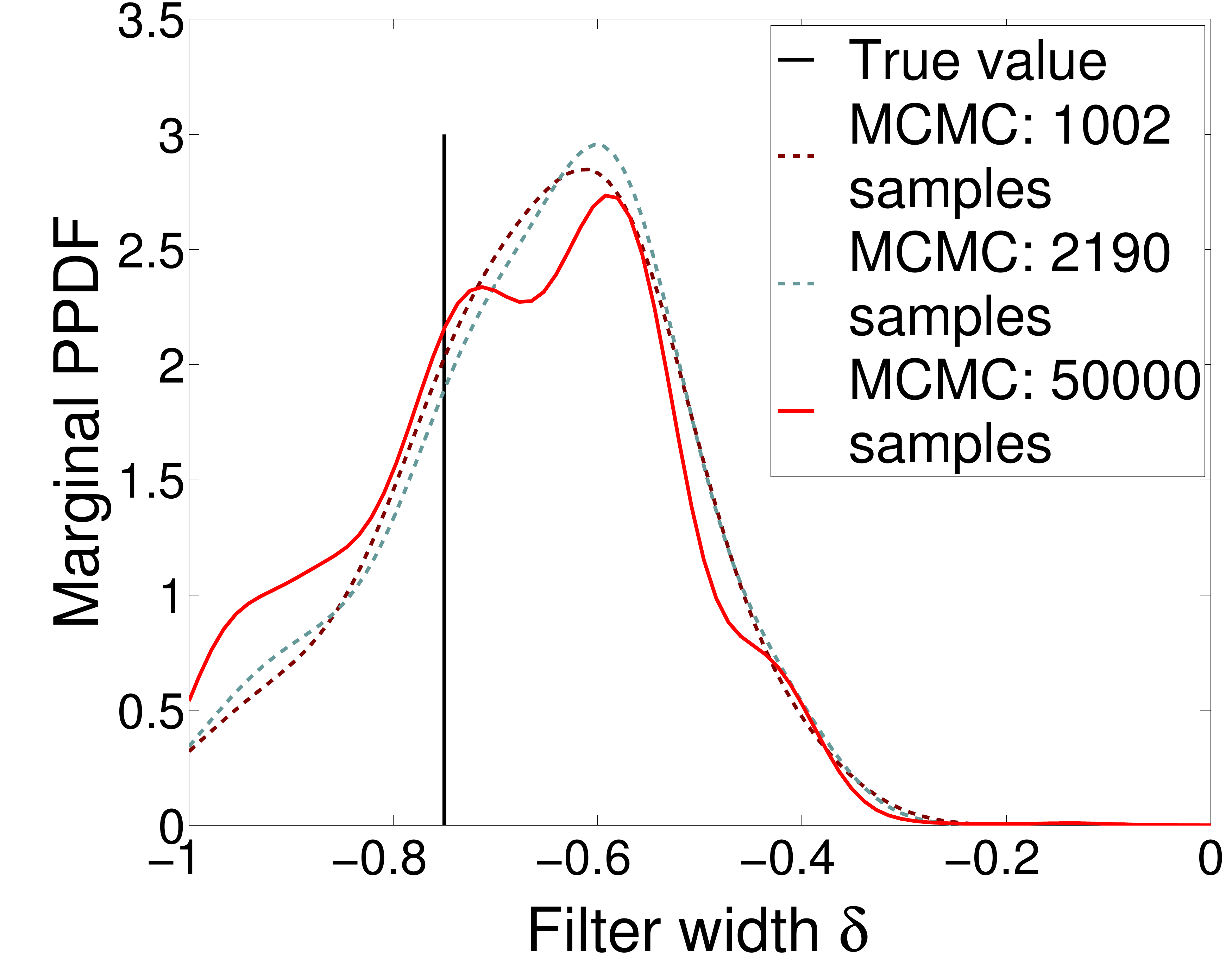} 
\caption{Marginal posterior probability density functions of model parameters with \eqref{EXP} likelihood function estimated using the linear surrogate systems on level 8 adaptive sparse grids with 1002, 2190 and 50000 samples {(excluding 10000 samples for burn-in period). These are the numbers of {model executions} that the conventional MCMC requires to obtain similar results.}}
\label{fig:PPDF4}
\vspace{-.1in}
\end{figure}

In order to demonstrate that our adaptive refinement strategy based on the smoothness of output data in probability space allows the change of likelihood models with minimal computational cost, we perform the above experiment with \eqref{EXP} likelihood function and $\xi = 500$ using the same surrogate of outputs. The marginal PPDFs of model parameters estimated using the linear surrogate systems are shown in Figure \ref{fig:PPDF3}. Again, they can be compared with marginal PPDFs estimated using conventional MCMC with the same number of model executions in Figure \ref{fig:PPDF4}. The plots confirm the accuracy of the surrogate PPDF for MCMC simulations and that surrogate-based MCMC requires less forward model executions than the conventional approach. On the other hand, it should be noted that some likelihood models, especially those resulting in peaky PPDFs, may require a surrogate system more accurate than that on level $8$ sparse grid. In those cases, the surrogate needs to be constructed on a grid of higher level. 

\begin{figure}[h!]
\centering
\hspace{0.6cm}\includegraphics[height=4.8cm]{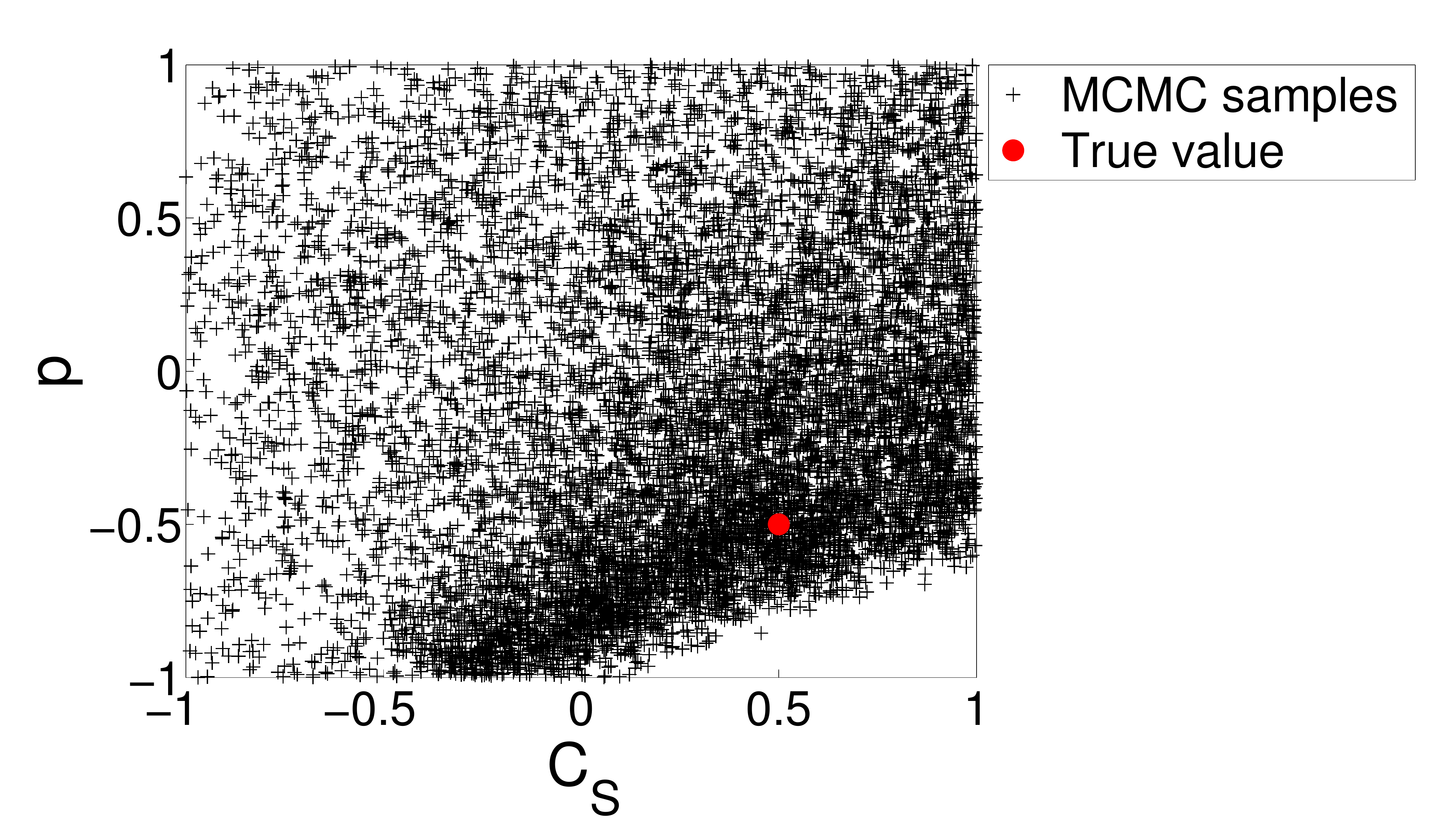} 
\caption{Two-dimensional marginal posterior probability density function of $C_S$ and $p$ with \eqref{MVN} likelihood model. The MCMC samples are obtained using the linear surrogate system on level $8$ sparse grid.}
\label{fig:scatter}
%\vspace{-.1in}
\end{figure}

The calibration results for both likelihood models show that the Smagorinsky constant $C_S$ and van Driest damping parameter $p$ have posterior maximizers near their true values, while smaller values are somewhat preferred for the filter width $\delta$. Meanwhile, the posterior distribution of $\delta$ is peaky, indicating that the data depend on $\delta$ and the Smagorinsky models with our selections of filter width (and spatial resolution) are incomplete. Indeed, finer grids are needed to sufficiently resolve the energy. The plots also reveal that the data are significantly more sensitive with respect to $\delta$ than to other parameters. This elucidates why finding the optimal value for $\delta$, i.e., determining the ideal place to truncate scale, is a very important issue in LES practice. Finally, the positive correlation between $C_S$ and $p$ can be observed in Figure \ref{fig:scatter}, in which the posterior samples projected on the $(C_S,p)$-plane are plotted. Given that our calibration data are extracted in near wake region, this correlation is expected. As larger value of $C_S$ increases the Smagorinsky lengthscale $\ell_S$, larger $p$ would be needed for a stronger damping of $\ell_S$ near the boundary. 

\section{Conclusion} 

\label{sec:conclusion}

In this paper, we present a surrogate modeling approach based on the AHSG method for Bayesian inference, with application to quantification of parametric uncertainty of LES turbulence models. The method {is based on those developed in [62] for less complex geophysical and groundwater models}, is model independent and can be flexibly used together with any MCMC algorithm and likelihood function. The accuracy and efficiency of our approach is illustrated by virtue of the numerical example consisting of the Smagorinsky model of two-dimensional flow around a cylinder. We combine the hierarchical linear basis and the local adaptive sparse grid technique to construct surrogate systems with a small number of model executions. Although the forward model investigated herein is highly nonlinear and more complicated than those in previous studies, our analysis indicates that the surrogate system is accurate for reasonable specifications of search regions. Compared to the conventional MCMC simulation, our surrogate-based approach requires significantly less model executions for estimating the parameter distribution and quantifying predictive uncertainty. Given the extremely high cost of turbulence simulations, this computational efficiency is critical for the feasibility of Bayesian inference in turbulence modeling. 

While the performance of surrogate modeling method is evaluated in this work for a synthetic cylinder flow model on relatively coarse grids, we expect comparable results for practical, more complicated calibration and prediction problems using real-world data; since three-dimensional, more refined simulations and real experiments of these flows are known to produce similar patterns to the 
investigated physical outputs in this study.
Still, a {three-dimensional} demonstration of our surrogate-based approach for these problems is irreplaceable and would be the next logical step. The framework presented here could be directly applied to other engineering flow models, as well as to the tasks of quantifying the structural uncertainties and comparing competing turbulence closure models. The accuracy of surrogate-based MCMC in these cases needs to be tested, but the verification, which is much less computational demanding than running the conventional MCMC, is possibly worthwhile. {Finally, besides our AHSG, several other methods can be employed to construct the surrogate system. A thorough comparative assessment with those methods is essential to fully justify the efficiency of our approach in turbulence uncertainty quantification problems and would be considered in the future.}        

Concerning sparse grid interpolation methods, additional research in accelerating the convergence rate of the surrogate is necessary. {One} direction is high-order sparse grid methods, which utilize high-order (instead of linear) hierarchical polynomial basis and whose superior efficiency has been justified for uncertainty quantification of groundwater models \cite{ZLY+13}. On the other hand, given that the outputs and PPDFs {do} not experience same level of sensitivity to different calibration parameters, combining {locally grid refinement strategy with dimension-adaptive sparse grid methods} to further reduce the number of interpolation points is {worth studying}.


\begin{thebibliography}{abbrv}

%\bibitem[AL99]{AL99} \textsc{N.A. Adams and A. Leonard}, \emph{Deconvolution
%of subgrid scales for the simulation of shock-turbulence interaction}, p.
%201 in: \textit{Direct and Large Eddy Simulation III}, (eds.: P. Voke, N.D.
%Sandham and L. Kleiser), Kluwer, Dordrecht, 1999.

%\bibitem[AL07]{AL07}\textsc{M. Anitescu, W. Layton}, \emph{Sensitivities in Large Eddy Simulation and Improved Estimates of Turbulent Flow Functionals},
%SIAM Journal on Scientific Computing 01/2007; 29:1650-1667. 

%\bibitem[ALP04]{ALP04} \textsc{M. Anitescu, W. Layton and F. Pahlevani}, 
%\emph{Implicit for local effects, explicit for nonlocal is unconditionally
%stable},ETNA, 18[2004] 174-187.

\bibitem{Ant81} \textsc{M. Antonopoulos-Domis}, \emph{Large eddy simulation of a passive scalar in isotropic turbulent}, Journal of Fluid Mechanics 104, 55--79, 1981.

\bibitem{BNT07}\textsc{ I. Babuska, F. Nobile, and R. Tempone}, \emph{A stochastic collocation method for elliptic partial differential equations with random input data}, SIAM J. Numer. Anal. 45 (2007), no. 3, 1005-1034 (electronic). MR MR2318799 (2008e:65372)

\bibitem{BNZ04}\textsc{I. Babuska, F. Nobile, and E. Zouraris}, \emph{Galerkin Finite Element Approximations of Stochastic Elliptic Partial Differential Equations}, SIAM Journal of Numerical Analysis, Vol. 42, pp. 800- 825, 2004.

\bibitem{BCWZ14} \textsc{F. Bao, Y. Cao, C.~G.~Webster and G. Zhang}, \emph{A hybrid sparse grid approach for nonlinear filtering problems based on adaptive domain approximations of the {Z}akai equation}, SIAM J. on Uncertainty Quantification, 2: 784-804, 2014.

\bibitem{BM94} \textsc{P. Beaudan, P. Moin}, \emph{Numerical Experiments on the Flow Past a Circular Cylinder at Sub-critical Reynolds Number}, in Report No. TF-62, Thermosciences Division, Department of Mechanical Engineering,  Stanford University, 1994.  

\bibitem{BW72} \textsc{E. Berger, R. Wille}, \emph{Periodic flow phenomena}, Ann. Rev. Fluid Mech. 4, 313-340, 1972.

\bibitem{BIL04} \textsc{L.C. Berselli, T. Iliescu and W. Layton}, \emph{Large Eddy Simulation}, Springer, Berlin, 2004.

%\bibitem[B61]{B61} \textsc{R. Betchov}, \emph{Semi-isotropic turbulence and
%helicoidal flows}, Phys. Fluids, 4(1961) 925-926.

\bibitem{BB92} \textsc{K. Beven and A. Binley}, \emph{The future of distributed models - Model calibration and uncertainty prediction}, Hydrol. Processes, 6(3), 279--298, 1992.

\bibitem{BR81} \textsc{S. Biringen, W.C Reynolds}, \emph{Large eddy simulation of the shear-free turbulent boundary layer}, Journal of Fluid Mechanics 103, 53--63, 1981.

%\bibitem[BIR09]{BIR09} \textsc{J. Borggaard, T. Iliescu and J.P. Roop}, 
%\emph{A Bounded Artificial Viscosity Large Eddy Simulation Model}, SIAM J.
%Numer. Anal., 47(2009), 622-645.

%\bibitem[BRTT12]{BRTT12} \textsc{A. Bowers, L. Rebholz, A. Takhirov, and C.
%Trenchea}, \emph{Improved accuracy in regularization models of
%incompressible flow via adaptive nonlinear filtering}, IJNMF, to appear,
%2012.

%\bibitem{Bou77} \textsc{J. Boussinesq}, \emph{Essai sur la th\'{e}orie des eaux courantes}, M\'{e}m. pr\'{e}s par div. savants \`{a} la Acad. Sci., 23, 1-680, 1877.

\bibitem{BT92} \textsc{G. Box and G. Tiao}, \emph{Bayesian Inference in Statistical Analysis}, 608
pp., Wiley-Interscience, N. Y., 1992.

%\bibitem[B98]{Boyd98} \textsc{J.P. Boyd}, \emph{Two comments on filtering
%for Chebyshev and Legendre spectral and spectral element methods :
%Preserving the boundary conditions and interpretation of the filter as adiffusion}, JCP, 143 (1998), 283-288.

\bibitem{Bre98} \textsc{M. Breuer}, \emph{Large eddy simulation of the subcritical flow past a circular cylinder: Numerical and modeling aspects}, Int. J. Numer. Meth. Fluids 28, 1281--1302, 1998.

\bibitem{BG04} \textsc{H.~J.~Bungartz, M.~Griebel},  \emph{{Sparse grids}}, Acta Numerica, 13, 1--123, 2004.


\bibitem{COP+11}\textsc{S.H. Cheung, T.A. Oliver, E.E. Prudencio, S. Prudhomme, R.D Moser} \emph{Bayesian inference with applications to turbulence modeling}, Reliab. Eng. Syst. Safety., 96, 1137--1149, 2011.

\bibitem{CFR79} \textsc{R.A. Clark, J.H. Ferziger, and W.C. Reynolds}, \emph{Evaluation of subgrid-scale models using an accurately simulated turbulent flow}, Journal of Fluid Mechanics 91, 1--16, 1979.

%\bibitem[CL10]{CL10} \textsc{J. Connors and W. Layton}, \emph{On the
%accuracy of the finite element method plus time relaxation}, Math. Comp., 79
%(2010), 619-648.

\bibitem{Dea70} \textsc{J.W. Deardorff}, \emph{A three-dimensional numerical study of turbulent channel flow at large Reynolds numbers}, Journal of Fluid Mechanics 41, 453--480, 1970.

\bibitem{DW11}\textsc{E. Dow and Q. Wang}, \emph{Quantification of structural uncertainties in the $k-\omega$ turbulence model}. AIAA Paper, 2011-1762, 2011.

%\bibitem[D04]{D04} \textsc{A. Dunca}, \emph{Space averaged Navier-Stokes
%equations in the presence of walls}, PhD Thesis, University of Pittsburgh,
%2004.

%\bibitem[D02]{D02} \textsc{A. Dunca,} \emph{Investigation of a shape
%optimization algorithm for turbulent flows}, report ANL/MCS-P1101-1003,
%Argonne National Lab, 2002, http://www-fp.mcs.anl.gov/division/publications/.

%\bibitem[DE04]{DE04} \textsc{A. Dunca and Y. Epshteyn}, \emph{On the
%Stolz-Adams de-convolution LES\ model}, SIAM J. Math. Anal., 37 (2006),
%1890-1902.

\bibitem{ELI13} \textsc{M. Emory, J. Larsson and G. Iaccarino}, \emph{Modeling of structural uncertainties in Reynolds-averaged Navier-Stokes closures}, Phys. Fluids 25, 110822, 2013.

%\bibitem[ELN07]{ELN07} \textsc{V. Ervin, W. Layton and M. Neda} ,\emph{\
%Numerical analysis of a higher order time relaxation model of fluids}, Int.
%J. Numer. Anal. and Modeling, 4 (2007), 648-670.

%\bibitem[ELN10]{ELN10} \textsc{V. Ervin, W. Layton and M. Neda},\emph{\
%Numerical analysis of filter based stabilization for evolution equations},
%technical report, TR-MATH 10-01, 2010,
%http://www.mathematics.pitt.edu/research/technical-reports.php.

%\bibitem[FM01]{FM01} \textsc{P. Fischer and J. Mullen}, \emph{Filter-based
%stabilization of spectral element methods}, C. R. Acad. Sci. Paris, 332(1),
%265 (2001).

\bibitem{FWK08} \textsc{J. Foo, X. Wan, and G. Karniadakis}, \emph{The multi-element probabilistic collocation method (ME-PCM): Error analysis and applications}, Journal of Computational Physics, 227, 9572--9595, 2008.

\bibitem{GL06} \textsc{D. Gamerman, and H. Lopes}, \emph{Markov Chain Monte Carlo : Stochastic Simulation for Bayesian Inference}, 2nd ed., 344 pp., Chapman and Hall, London, 2006.

\bibitem{GZ07} \textsc{B. Ganapathysubramanian, N. Zabaras}, \emph{Sparse grid collocation schemes for stochastic natural convection problems}, J. Comput. Phys. 225(1): 652-685, 2007.

\bibitem{GG03} \textsc{T. Gerstner, and M. Griebel}, \emph{Dimension-Adaptive Tensor-Product Quadrature}, Computing, 71 (2003), pp. 65-87.

\bibitem{GS91} \textsc{R. G. Ghanem and P. D. Spanos}, \emph{Stochastic Finite Elements: A Spectral Approach}, Springer-Verlag, New York, NY, 1991.

\bibitem{GELI12} \textsc{C. Gorl\'{e}, M. Emory, J. Larsson and G. Iaccarino}, \emph{Epistemic uncertainty quantification for RANS modeling of the flow over a wavy wall}, Center for Turbulence Research, Annual Research Briefs 2012.

\bibitem{GI13} \textsc{C. Gorl\'{e} and G. Iaccarino}, \emph{A framework for epistemic uncertainty quantification of turbulent scalar flux models for Reynolds-averaged Navier-Stokes simulations}, Phys. Fluids 25, 055105, 2013.

\bibitem{Gri98} \textsc{M. Griebel}, \emph{Adaptive sparse grid multilevel methods for elliptic PDEs based on finite differences}, Computing, 61(2), 151--179, 1998. doi:10.1007/BF02684411.

%\bibitem{Gri12} \textsc{M. Grigoriu}, \emph{A Method for Solving Stochastic Equations by Reduced Order Models and Local Approximations}, Journal of Computational Physics, Vol. 231, pp. 6495-6513, 2012.

%\bibitem[GAS09]{GAS09} \textsc{E. Garnier, N. Adams and P. Sagaut},\emph{\
%Large eddy simulation for compressible flows}, Springer, Berlin, 2009.

\bibitem{GWZ13b} \textsc{M. Gunzburger, C.~G.~Webster and G. Zhang}, \emph{An adaptive sparse grid iterative ensemble {K}alman filter approach for parameter field estimation}, Inter. J. of Comp. Math., 91(4): 798-817, 2014.

\bibitem{GWZ14} \textsc{M. Gunzburger, C.~G.~Webster and G. Zhang}, \emph{Stochastic finite element methods for partial differential equations with random input data}, Acta Numerica, 23: 521--650, 2014.

\bibitem{HLMS06} \textsc{H. Haario, M. Laine, A. Mira, and E. Saksman}, \emph{DRAM: Efficient adaptive MCMC}, Stat. Comput., 16(4), 339--354, 2006. doi:10.1007/s11222- 006-9438-0.

\bibitem{Hec12} \textsc{F. Hecht}, \emph{New development in freefem++}. J. Numer. Math. 20, 251--265, 2012.

%\bibitem[HePi]{HePi} \textsc{F. Hecht and O. Pironneau}, \emph{FreeFEM++ },
%webpage: http://www.freefem.org.

%\bibitem[HWM88]{HWM88} \textsc{J.C. Hunt, A.A. Wray and P. Moin}, \emph{%
%Eddies stream and convergence zones in turbulent flows}, CTR report CTR-S88,
%1988.

\bibitem{John04} \textsc{V. John}, \emph{Reference values for drag and lift of a two-dimensional time-dependent flow around a cylinder}, Int. J. Numer. Math. Fluids, 44, 777--788, 2004.

\bibitem{JM01} \textsc{V. John, G. Matthies}, \emph{Higher-order finite element discretizations in a benchmark problem for incompressible flows}, Int. J. Numer. Meth. Fluids, 37, 885--903, 2001.

%\bibitem{KW05} \textsc{ A. Klimke and B. Wohlmuth}, \emph{Algorithm 847: Spinterp: piecewise multilinear hierarchical sparse grid interpolation in MATLAB}, ACM Transactions on Mathematical Software (TOMS) 2005, 31(4), 561--579.

\bibitem{Lilly67} \textsc{D. K. Lilly}, \emph{The representation of small scale turbulence in numerical simulation experiments}. In H.H. Goldstine, editor, Proc. IBM Sci. Computing Symp. On Environmental Sciences, pages 195-210, Yorktown Heights, NY, 1967.

\bibitem{LCK10} \textsc{X. Liu, M. A. Cardiff, and P. K. Kitanidis}, \emph{Parameter estimation in nonlinear environmental problems}, Stochastic. Environ. Res. Risk Assess., 24(7), 1003--1022, 2010.

\bibitem{LGH11} \textsc{M. Liu, Z. Gao, and J. Hesthaven}, \emph{Adaptive Sparse Grid Algorithms with Applications to Electromagnetic Scattering under Uncertainty}, Applied Numerical Mathematics, Vol. 61, pp. 24-37, 2011.

%\bibitem[LMNR08]{LMNR08} \textsc{W. Layton, C Manica, M Neda, L Rebholz}, 
%\emph{Numerical analysis and computational testing of a high accuracy
%Leray-Deconvolution model of turbulence}. NMPDEs 2008; 24(2):555--582.

%\bibitem[L07b]{L07b} \textsc{W. Layton} , \emph{Superconvergence of finite
%element discretization of time relaxation models of advection}, BIT, 47
%(2007), 565-576.

%\bibitem[LMNR06]{LMNR06} \textsc{W. Layton, C. Manica, M. Neda and L. Rebholz%
%}, \emph{The joint Helicity-Energy cascade for homogeneous, isotropic
%turbulence generated by approximate deconvolution models,} Adv. and Appls.
%in Fluid Mechanics, 4 (2008), 1-46.

%\bibitem[LN07]{LN07} \textsc{W. Layton and M. Neda}, \emph{Truncation of
%scales by time relaxation}, JMAA, 325 (2007), 788-807.

%\bibitem[Loe10]{Loe10} \textsc{G. J. A. Loeven}, \emph{Efficient uncertainty quantification in computational fluid dynamics}, PhD Thesis, Technische Universiteit Delft, 2010.

%\bibitem[LT83]{LT83} \textsc{E. Levich and A. Tsinober}, \emph{On the role
%of helical structures in 3-dimensional turbulent flows}, Physics Letters 93A
%(1983) 293-297.

%\bibitem[LT83b]{LT83b} \textsc{E. Levich and A. Tsinober}, \emph{Helical
%structures, fractal dimensions and renormalization group approach in
%homogeneous turbulence}, Physics Letters 96A (1983) 292-297.

\bibitem{LMS07}\textsc{D. Lucor, J. Meyers, P. Sagaut}, \emph{Sensitivity analysis of large-eddy simulations to subgrid-scale-model parametric uncertainty using polynomial chaos}, J. Fluid Mech. 585, 255-279, 2007.

%\bibitem[Lilly86]{Lilly86} \textsc{D.K. Lilly}, \emph{The structure,
%energetics and propagation of rotating convection storms, II: Helicity and
%Storm Stabilization}, J Atmos. Sciences, 43(1986),126-140.

\bibitem{MZ09} \textsc{X. Ma, and N. Zabaras}, \emph{An efficient Bayesian inference approach to inverse problems based on an adaptive sparse grid collocation method}, Inverse Probl., 25(3), 035013, 2009.

%\bibitem{MZ10} \textsc{X. Ma, and N. Zabaras}, \emph{An adaptive high-dimensional stochastic model representation technique for the solution of stochastic partial differential equations}, Journal of Computational Physics 229 (2010) 3884--3915.

\bibitem{Mas89} \textsc{P.J. Mason}, \emph{Large eddy simulation of the convective atmospheric boundary layer}, Journal of the Atmospheric Sciences 46, 1492--1516, 1989.

\bibitem{MD90} \textsc{P.J. Mason, S.H. Derbyshire}, \emph{Large-eddy simulation of the stable-stratified atmospheric boundary layer}, Boundary-layer Meteorology 53, 117-162, 1990.

\bibitem{MK82} \textsc{P. Moin, J. Kim}, \emph{Numerical investigation of turbulent channel flow}, Journal of Fluid
Mechanics 18, 341--377, 1982.

\bibitem{MMRF79} \textsc{N.N. Monsour, P. Moin, W.C. Reynolds, and J.H. Ferziger}, \emph{Improved methods for large eddy simulations of turbulence}, In Turbulent Shear Flows I, 386--401, 1979.

%\bibitem[MLF03]{MLF03} \textsc{J. Mathew, R. Lechner, H. Foysi, J.
%Sesterhenn and R. Friedrich}, \emph{An explicit filtering method for large
%eddy simulation of compressible flows}, Physics of Fluids, 15 (2003),
%2279-2289.

%\bibitem[MF98]{MF98} \textsc{J.S. Mullen and P.F. Fischer}, \emph{Filtering
%techniques for complex geometry fluid flows}, Comm. in Num. Meth. in Eng.,
%15 (1999), 9-18.

%\bibitem[ND99]{WALE99} \textsc{F. Nicoud and F. Ducros}, \emph{Subgrid-Scale
%Stress Modelling Based on the Square of the Velocity Gradient Tensor}, Flow,
%Turbulence and Combustion, 62 (1999), 183-200.

\bibitem{NR50} \textsc{J. von Neumann and R.D. Richtmyer}, \emph{A method for the numerical calculation
of hydrodynamic shocks}. J. Appl. Phys., 21, 232-237, 1950.

\bibitem{NTW08} \textsc{F. Nobile, R. Tempone, C. G. Webster}, \emph{A Sparse Grid Stochastic Collocation Method for Partial Differential Equations with Random Input Data}, SIAM J. Numer. Anal., 46(5), 2309--2345, 2008.

\bibitem{NTW08b} \textsc{F. Nobile, R. Tempone, and C. G. Webster}, \emph{An anisotropic sparse grid collocation method for elliptic partial differential equations with random input data}, SIAM J. Numer. Anal. 46 (5), 2411--2442,  2008.

\bibitem{Nor87} \textsc{C. Norberg}, \emph{Effects of Reynolds number and a low-intensity free-stream turbulence on the flow around a circular cylinder}, Publication No. 87/2, Department of Applied Thermodynamics and Fluid Mechanics, Chalmer University of Technology, Gothenburg, Sweden.

\bibitem{OM11}\textsc{T. Oliver and R. Moser}, \emph{Bayesian uncertainty quantification applied to RANS turbulence models}, Journal of Physics: Conference Series 318, 042032, 2011.

\bibitem{Pfl10} \textsc{D. Pfl\"{u}ger}, \emph{Spatially Adaptive Sparse Grids for High-Dimensional Problems}, Ph.D. Thesis, TU Munich, Munich, Germany, 2010.

\bibitem{Pope00} \textsc{S. Pope}, \emph{Turbulent flows}, Cambridge University Press, 2000.

\bibitem{Rah08} \textsc{S. Rahman}, \emph{A Polynomial Dimensional Decomposition for Stochastic Computing}, International Journal for Numerical Methods in Engineering, Vol. 76, 2008, pp. 2091--2116.

\bibitem{RTB12} \textsc{S. Razavi, B. A. Tolson, and D. H. Burn}, \emph{Review of surrogate modeling in water resources}, Water Resour. Res., 48, W07401, 2012.

\bibitem{RC04} \textsc{C. Robert, G. Casella}, \emph{Monte Carlo Statistical Methods}, 2nd ed., Springer, 2004. 

\bibitem{RFBP97} \textsc{W. Rodi, J.H. Ferziger, M. Breuer and M. Pourqui\'{e}}, \emph{Status of large eddy simulation: results of a workshop}, Workshop on LES of Flows Past Bluff Bodies, Rottach-Egern, Tegernsee, Germany, June 26-28, 1995, J. Fluids Eng., 119, 248--262, 1997.

%\bibitem[S01]{Sagaut} \textsc{P. Sagaut, }\emph{Large eddy simulation for
%Incompressible flows,} Springer, Berlin, 2001.

%\bibitem[SL97]{SL97} \textsc{P. Sagaut and T. Le}, \emph{Some investigations of the sensitivity of large eddy simulation}, Tech.
%Rep. 1997-12, ONERA, Toulouse, France, 1997.

%\bibitem[SV10]{SV10} \textsc{G. Schoups, and J. A. Vrugt}, \emph{A formal likelihood function for parameter and predictive inference of hydrologic models with correlated, heteroscedastic, and non-Gaussian errors}, Water Resour. Res., 46, W10531, 2010. doi:10.1029/2009WR008933.

\bibitem{SM11} \textsc{S. Sankaran, A. Marsden}, \emph{A stochastic collocation method for uncertainty quantification and propagation in cardiovascular simulations}, J. Biomech Eng 133(3), 031001, 2011.

%\bibitem[ST96]{ST96}\textsc{M. Sh\"{a}fer and S. Turek}, \emph{Benchmark computations of laminar flow around cylinder}, Flow Simulation with High-Performance Computers II, Vieweg, 1996.

\bibitem{Sma63}\textsc{J.S. Smagorinsky}, \emph{General circulation experiments with the primitive equations}, Mon. Weather Review, 91, 99-164, 1963.

\bibitem{SBT08} \textsc{P. Smith, K. J. Beven, and J. A. Tawn}, \emph{Informal likelihood measures in model assessment: Theoretic development and investigation}, Adv. Water Resour., 31, 1087--1100, 2008.

%\bibitem[SSM+10]{SSM+10} \textsc{T. Smith, A. Sharma, L. Marshall, R. Mehrotra, and S. Sisson}, \emph{Development of a formal likelihood function for improved Bayesian inference of ephemeral catchments}, Water Resour. Res., 46, W12 551, 2010. doi:10.1029/2010WR009514.

\bibitem{Sto13} \textsc{M. Stoyanov}, \emph{User Manual: TASMANIAN sparse grid}, ORNL Technical Report, 2013. 

\bibitem{TTW14}\textsc{H. Tran, C. Trenchea, C. Webster}, \emph{A convergence analysis of stochastic collocation method for Navier-Stokes equations with random input data}, ORNL Technical Report, 2014.  

\bibitem{vD56} \textsc{E.R. van Driest}, \emph{On turbulent flow near a wall}, J. Aerospace Sci., 23, 1007--1011, 1956.

\bibitem{VBC+08} \textsc{J.A. Vrugt, C. J. F. ter Braak, M. P. Clark, J. M. Hyman, and B. A. Robinson}, \emph{Treatment of input uncertainty in hydrologic modeling: Doing hydrology backward with Markov chain Monte Carlo simulation}, Water Resour. Res., 44, W00B09, 2008.

\bibitem{VBD+09}\textsc{J. Vrugt, C. Ter Braak, C. Diks, D. Higdon, B. Robinson, and J. Hyman}, \emph{Accelerating Markov chain Monte Carlo simulation by differential evolution with self-adaptive randomized subspace sampling}, Int. J. Nonlinear Sci. Numer. Simulation, 10, 273--290, 2009.

%\bibitem[WG13]{WG13} \textsc{Q. Wang and J.-H. Gao}, \emph{The drag-adjoint field of a circular cylinder wake at Reynolds numbers 20, 100 and 500},
% submitted, 2013. 
 
\bibitem{WK06} \textsc{X. Wan and G. E. Karniadakis}, \emph{Long-term behavior of polynomial chaos in stochastic flow simulations}, Comput. Methods Appl. Mech. Engrg. 195, 5582--5596, 2006.

\bibitem{WLSB08} \textsc{J.A.S. Witteveen, G.J.A. Loeven, S. Sarkar, and H. Bijl}, \emph{Probabilistic collocation for period-1 limit cycle oscillations}, Journal of Sound and Vibration 311(1-2), 421-439, 2008.

%\bibitem[ZSZW12]{ZSZW12} \textsc{L. Zeng, L. Shi, D. Zhang, and L. Wu}, \emph{A sparse grid based Bayesian method for contaminant source identification}, Adv. Water Resour., 37, 1--9, doi:10.1016/j.advwatres.2011.09.011.

%\bibitem[ZG12]{ZG12} \textsc{G. Zhang and M. Gunzburger}, \emph{Error analysis of a stochastic collocation method for parabolic
%partial differential equations with random input data}, SIAM J. Numer. Anal. 50 (2012), no. 4, 1922-1940.

\bibitem{ZLY+13} \textsc{G. Zhang, D. Lu, M. Ye, M. Gunzburger, C. G. Webster}, \emph{An adaptive sparse grid high-order stochastic collocation method for Bayesian inference in ground water reactive transport modeling},  Water Resources Research, 49, 6871--6892, 2013.
\end{thebibliography}
\end{document}